\documentclass[twocolumn]{autart}

\usepackage{graphicx}
\usepackage{amssymb}
\usepackage{amsmath}
\usepackage{bm}
\usepackage{tikz}

\usepackage{enumitem}

\usepackage{algorithm}
\usepackage{algorithmicx}
\usepackage{algpseudocode}

\usetikzlibrary{external}
\usepackage[round]{natbib}

\usepackage[hidelinks]{hyperref}
\hypersetup{
     colorlinks   = true,
     citecolor    = blue 
}
\begin{document}
    \begin{frontmatter}
        \title{Numerical methods for solving minimum-time problem for linear systems\thanksref{footnoteinfo}} 
        
        \thanks[footnoteinfo]{This paper was not presented at any IFAC 
        meeting. Corresponding author M.~E.~Buzikov. Tel. +7 (966) 095-03-04.}
        
        \author[ICS]{Maksim E. Buzikov}\ead{me.buzikov@physics.msu.ru},
        \author[ICS]{Alina M. Mayer}\ead{mayer@ipu.ru}
        
        \address[ICS]{V.~A.~Trapeznikov Institute of Control Sciences of Russian Academy of Sciences, Moscow, Russia} 
                  
        
        \begin{abstract}
            This paper offers a contemporary and comprehensive perspective on the classical algorithms utilized for the solution of minimum-time problem for linear systems (MTPLS). The use of unified notations supported by visual geometric representations serves to highlight the differences between the Neustadt-Eaton and Barr-Gilbert algorithms. Furthermore, these notations assist in the interpretation of the distance-finding algorithms utilized in the Barr-Gilbert algorithm. Additionally, we present a novel algorithm for solving MTPLS and provide a constructive proof of its convergence. Similar to the Barr-Gilbert algorithm, the novel algorithm employs distance search algorithms. The design of the novel algorithm is oriented towards solving such MTPLS for which the analytic description of the reachable set is available. To illustrate the advantages of the novel algorithm, we utilize the isotropic rocket benchmark. Numerical experiments demonstrate that, for high-precision computations, the novel algorithm outperforms others by factors of tens or hundreds and exhibits the lowest failure rate.
        \end{abstract}
    \end{frontmatter}

    \section{Introduction}

    Since the discovery of the maximum principle, the field of optimal control has experienced considerable and accelerated progress, particularly in the development of numerical methods. The maximum principle provides a foundational framework that enables researchers to systematically approach and solve a wide range of optimal control problems. Among these, the application of numerical techniques to address minimum-time problems for linear systems (MTPLS) has proven successful. For such problems, algorithms for calculating the optimal control parameters have been developed, and the global convergence of some of these algorithms has been proven. Furthermore, several studies have been conducted to increase the convergence rate of these algorithms.

    Currently, two successful methods are available for developing algorithms that calculate the solution to the MTPLS\footnote{The focus of our attention will be directed only to problems with control input constraints of the inclusion type. Integral control input constraints are the topic of a separate study}. Both methods are more readily understandable when geometric interpretations of the reachable set are employed. The reachable set is defined as the collection of all potential states that a dynamic system can attain from a given initial state through the application of a set of admissible control inputs over a specified time horizon. In linear systems, the reachable set is convex, and it is relatively straightforward to obtain a parameterization of its boundary using a support vector.

    The first method was initially proposed by \citet{Neustadt1960-ki} and entails reducing the MTPLS to a finite-dimensional maximization of a specific function. In the following, this function will be referred to as the \textit{boosting-time function}. This function takes a support vector as input and returns the minimal moment of time when the target state hits the tangent hyperplane to the reachable set with the given support vector. \citet{Neustadt1960-ki} demonstrated that the maximum of the boosting-time function is equal to the value of the MTPLS. Furthermore, he showed that the support vector which maximizes the function defines the optimal control. \citet{Eaton1962-ai} and later \citet{Neustadt1963-en} described a steepest ascent algorithm for maximizing the boosting-time function. From a geometric perspective, the selection of the step size for the steepest ascent algorithm is analogous to the inclination of the tangent hyperplane that separates the target state from the reachable set. With each iteration of the approximation of the support vector, a new value for the boosting-time function is established, forming an increasing sequence. By alternating the inclination of the tangent hyperplane and the time increment, sequences of support vectors and time moments are obtained that are expected to approach the optimal values. In order to guarantee the convergence of the algorithm, \citet{Boltyanskii1971-ty} provided constructively verifiable conditions for checking the step size (degree of inclination of the tangent hyperplane). The results presented by \citet{Neustadt1960-ki} and \citet{Eaton1962-ai} were independently rediscovered by \citet{Pshenichnyi1964-ux}. The generalized approach was further developed by \citet{Kirin1964-vn}, \citet{Babunashvili1964-tp}, \cite{Ivanov1971-pe}, and \cite{Barnes1972-oc,Barnes1973-ie}.

    \citet{Paiewonsky1963-kg} observed that the steepest ascent method may demonstrate a slow convergence in certain instances. They suggested employing Powell's method to enhance the convergence rate. \citet{Fadden1964-rk} outlined the issues with convergence and suggested ad-hoc solutions. \citet{Pshenichnyi1968-ua} adapted the method of \citet{Fletcher1963-rp} to solve the MTPLS. Additionally, they proposed approximating the boosting-time function by a quadratic polynomial to determine the step size of the algorithm. The final algorithm of \cite{Pshenichnyi1968-ua} is quite complex due to the presence of backtracks in exceptional cases, yet its convergence rate is superior to that of previous algorithms. It is also noteworthy that the final version of their algorithm lacks a consistent proof of convergence, as the algorithm relies on approximations.

    The second method of solving MTPLS was proposed by \citet{Ho1962-ki} and consists of calculating the distance from the target state to the reachable set and varying the time to identify the minimum moment when the distance becomes zero. The majority of works devoted to this method focus on the problem of computing the distance from the target state to the reachable set. The nature of MTPLS allows the boundary of the reachable set to be parameterized by utilizing a \textit{contact function}. The contact function serves to parameterize the boundary of the convex set by means of a support vector. A hyperplane defined by a support vector is a tangent hyperplane that touches the convex set at the contact set, which is given by the contact function for the aforementioned support vector. The contact function for the reachable set of a linear system can be readily derived from the maximum principle by integrating the state equations with the extremal control inputs. In a notable contribution, \citet{Ho1962-ki} proposed an iterative procedure for computing the distance between the reachable set and the target state. \citet{Fancher1965-xk} modified this procedure to facilitate implementation on a computer. In order to demonstrate the convergence of the algorithm proposed by \citet{Fancher1965-xk}, \citet{Gilbert1966-sc} investigated the more general problem of determining the distance between a target state and an arbitrary convex set defined by a contact function. We will henceforth refer to this as the minimal distance problem, or MDP. As demonstrated by \citet{Gilbert1966-sc}, the minimal-distance algorithm proposed by \citet{Ho1962-ki} is based on a geometrical principle. In each iteration, the algorithm considers a line segment defined by the current approximation of the nearest point to the target state and the corresponding contact point. The nearest point to the target state on the line segment serves to define the subsequent approximation. \citet{Barr1966-tc,Barr1969-ca} proposed the use of a new approximation, namely a point from a convex polytope whose vertices are the previously calculated points of the reachable set. \citet{Barr1966-tc,Barr1969-ca} demonstrated through calculations that this approach markedly accelerates the convergence of the distance algorithm. The solution of the MDP defines the support vector such that the corresponding value of the contact function is the nearest point to the target state. To solve the MTPLS, \citet{Fujisawa1967-aj} and, independently, \citet{Barr1969-aj} proposed increasing the time variable using the boosting-time function with the support vector that is produced after calculation of the solution for the MDP.

    \citet{Demyanov1964-zf} independently proposed the same approach as \citet{Ho1962-ki} for solving the MTPLS and MDP. \citet{Baranov1965-pe} complemented these studies by examining the rate of convergence of the method for solving the MDP. In a subsequent contribution, \citet{Gabasov1966-bf} proposed the use of a boosting-time function in conjunction with Demyanov's algorithm to address the MTPLS. A more general formulation of the MDP was studied by \citet{Gindes1966-mw} and \citet{Demyanov1968-zo}, and a qualitative comparison of the aforementioned methods was made by \citet{Rabinovich1966-eg}.

    In the second method, the efficacy of solving the MTPLS depends on the speed of solving the MDP. The utilization of a convex polytope generated from previously computed points proved to be an exceptionally fruitful approach. In their seminal works, \citet{Barr1966-tc,Barr1969-ca} and \citet{Pecsvaradi1970-ro} proposed a set of guidelines for the selection of vertices for the convex polytope, which should be retained for the subsequent iteration. Moreover, it was demonstrated that these rules result in a reduction in the number of iterations. Perhaps the most successful selection rule was proposed by \citet{Gilbert1988-ug}. This rule entails retaining exclusively those vertices of the convex polytope that correspond to the face containing the point that is most proximate to the target state. This selection rule, when combined with the algorithm defined by \citet{Gilbert1966-sc}, constitutes an algorithm for solving the MDP, which is called the Gilbert–Johnson–Keerthi (GJK) distance algorithm. A generalization of the algorithm for arbitrary convex sets is provided by \citet{Gilbert1990-wx}. The effective and robust implementations of GJK proposed by \citet{Van-den-Bergen1999-fi} and \citet{Montanari2017-bw}. \citet{Qin2019-kw} employed a technique initially proposed by \citet{Nesterov1983-yw} to accelerate the \citet{Gilbert1966-sc} algorithm. Similarly, \citet{Montaut2024-dl} proposed a procedure to accelerate the GJK distance algorithm. \citet{Montaut2024-dl} also demonstrated that the GJK distance algorithm exhibits competitive performance in long-distance scenarios and remains a state-of-the-art algorithm.

    It should be noted that there are alternative methods for solving the MTPLS problem that do not rely on the \citet{Neustadt1960-ki} and \citet{Ho1962-ki} methods. An alternative iterative procedure, applicable to a more restricted class of problems for linear systems, was proposed by \citet{Knudsen1964-jt}. The application of this procedure requires the calculation of switching times, which is further complicated by the lack of constructively specified convergence conditions. Furthermore, the procedure is unable to be scaled to a more expansive class of linear problems. The same can be said about the methods proposed by \citet{Belolipetskii1977-he}, \citet{Gnevko1986-gi}, \citet{Kiselev1991-mq}, \citet{Aleksandrov2012-mw}.

    The objective of this paper is to present the methods and algorithms for solving MTPLS in a systematic manner, within a unified terminology and a general problem statement. This paper provides a general overview of algorithms based on Neustadt's method and Ho's method, which are capable of solving MTPLS with a moving target set. The primary contribution and innovation of this paper is the proposed technique for incorporating the potential for analytical calculation of the contact function for specific systems. We present a novel algorithm for solving MTPLS, based on Ho's method. Given the necessity of solving MDP within Ho's method, we have also developed a new algorithm for solving MDP that is comparable to the GJK distance algorithm but more straightforward to implement.

    The paper presents four algorithms for solving MDP and three algorithms for solving MTPLS. In order to facilitate a comparative analysis of the algorithms' performance, we will utilize the motion model proposed by \citet{Isaacs1955-bn} and commonly referred to as the "isotropic rocket". The isotropic rocket is a convenient benchmark for several reasons. First, the analytic description of the contact function for the reachable set of the isotropic rocket is known \citep{Buzikov2024-zm}. Second, the reachable set of the isotropic rocket is strictly convex but not strongly convex, which presents a challenge for the convergence rate of the algorithms.

    It is worth mentioning, that the solution we have derived to the problem for the isotropic rocket holds intrinsic value in its own right. The works of \citet{Bakolas2014-ze}, \citet{Buzikov2024-zm} present algorithms for solving the problem at an arbitrary final velocity of the isotropic rocket. \citet{Akulenko2011-kr} provided a description of a system of equations that can be solved to obtain parameters for optimal control in a problem with a fixed final velocity and a time-invariant target set. Nevertheless, it is challenging to propose methods for solving this system that will guarantee convergence to the optimal solution. Our algorithm is capable of identifying a guaranteed solution to the problem, even when the target set is time-varying. Another widely used model for fixing the velocity within the target set is the \citet{Dubins1957-hl} model. In the case of a stationary target set, an analytical solution can be obtained within the Dubins model. The issue of a moving target set has been addressed by \citet{Bakolas2013-zs,Gopalan2017-wy,Manyam2020-ym,Zheng2021-wu,Manyam2022-gh}. However, it cannot be stated that a convergent algorithm for arbitrary continuous motions of the target set has been developed. The isotropic rocket provides an alternative way to construct a reference trajectory for such cases.

    The paper is structured as follows. Section~\ref{sec:prel} provides preliminary information about convex sets and contact functions, which is essential to understand in order to contextualize the subsequent sections. Section~\ref{sec:prob_st} presents a general statement of the problem and outlines the assumptions that are necessary to justify the algorithms. Section~\ref{sec:constr} provides a general geometric overview of the existing results on the MTPLS. Additionally, Section~\ref{sec:constr} describes the essential building blocks required for the development of the algorithms. Section~\ref{sec:MDP} presents four algorithms for solving MDPs and proves the convergence of the novel algorithm for MDPs. Section~\ref{sec:MTPLS} describes three algorithms for solving the MTPLS and proves the convergence of the novel algorithm for the MTPLS. Section~\ref{sec:exp} outlines the practical considerations associated with the implementation of the proposed algorithms. Additionally, it describes the "isotropic rocket" benchmark and compares the performance of the algorithms. Finally, Section~\ref{sec:concl} concludes the paper with a discussion and brief description of future work.

    \section{Preliminaries}\label{sec:prel}

    We describe the state space by the column vectors with $n \in \mathbb{N}$ elements, i.e. the state space is $\mathcal{S} = \mathbb{R}^{n \times 1}$. The inner product and the norm are determined by:
    \begin{equation*}
        (\mathbf{a}, \mathbf{b}) \overset{\mathrm{def}}{=} \mathbf{a}^\top\mathbf{b}, \quad \lVert\mathbf{a}\rVert \overset{\mathrm{def}}{=} \sqrt{(\mathbf{a}, \mathbf{a})}, \quad \mathbf{a}, \mathbf{b} \in \mathcal{S}.
    \end{equation*}
    In the conjugate space $\mathcal{S}^* = \mathbb{R}^{1 \times n}$, we define
    \begin{equation*}
        (\mathbf{a}, \mathbf{b}) \overset{\mathrm{def}}{=} \mathbf{a}\mathbf{b}^\top, \quad \lVert\mathbf{a}\rVert \overset{\mathrm{def}}{=} \sqrt{(\mathbf{a}, \mathbf{a})}, \quad \mathbf{a}, \mathbf{b} \in \mathcal{S}^*.
    \end{equation*}

\begin{figure}
        \begin{center}
            \includegraphics[width=0.30\textwidth]{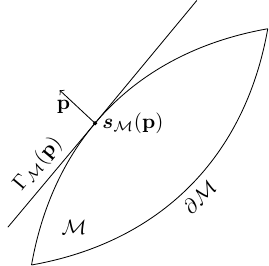}
            \caption{The support vector $\mathbf{p}$, the contact point $\boldsymbol{s}_{\mathcal{M}}(\mathbf{p})$, and the tangent hyperplane $\Gamma_{\mathcal{M}}(\mathbf{p})$}
            \label{fig:strictly_conv_comp}
        \end{center}
    \end{figure}

    In the further description we actively use properties of strictly convex compact sets. If $\mathcal{M} \subset \mathcal{S}$ is a strictly convex compact set then there exists a function $\boldsymbol{s}_{\mathcal{M}}: \mathcal{S}^* \to \mathcal{S}$ such that $\boldsymbol{s}_{\mathcal{M}}(\mathbf{p})$ is a {\it contact point} of $\mathcal{M}$ with the {\it support vector} $\mathbf{p} \in \mathcal{S}^*$, $\mathbf{p} \neq \boldsymbol{0}$, i.e. the hyperplane
    \begin{equation*}
        \Gamma_{\mathcal{M}}(\mathbf{p}) \overset{\mathrm{def}}{=} \{\mathbf{s} \in \mathcal{S}:\: \mathbf{p}(\mathbf{s} - \boldsymbol{s}_{\mathcal{M}}(\mathbf{p})) = 0\}
    \end{equation*}
    is a {\it tangent hyperplane} to $\mathcal{M}$ at $\boldsymbol{s}_{\mathcal{M}}(\mathbf{p})$ (Fig.~\ref{fig:strictly_conv_comp}). The function $\boldsymbol{s}_{\mathcal{M}}$ is a {\it contact function}.  In general, the contact function is given by
    \begin{equation}\label{eq:cont_def}
        \boldsymbol{s}_{\mathcal{M}}(\mathbf{p}) \overset{\mathrm{def}}{=} \mathrm{arg}\max_{\mathrm{s} \in \mathcal{M}}(\mathbf{p}\mathbf{s}).
    \end{equation}
    The $\mathrm{arg}\max$ exists due to compactness of $\mathcal{M}$. Its value is unique because $\mathcal{M}$ is strictly convex. The boundary of $\mathcal{M}$ can be expressed by the following:
    \begin{equation}\label{eq:dM}
        \partial\mathcal{M} = \{\boldsymbol{s}_{\mathcal{M}}(\mathbf{p}):\: \mathbf{p} \in \mathcal{S}^*,\: \mathbf{p} \neq \boldsymbol{0}\}.
    \end{equation}
    This representation is a parametrization of the boundary of the convex set by a support vector.
    
    The general property of the strictly convex set $\mathcal{M}$ is that the inequality
    \begin{equation}\label{eq:conv}
        \mathbf{p}(\boldsymbol{s}_{\mathcal{M}}(\mathbf{p}) - \mathbf{s}) \geq 0
    \end{equation}
    holds for any $\mathbf{p} \in \mathcal{S}^*$ ($\mathbf{p} \neq 0$) and $\mathbf{s} \in \mathcal{M}$. Moreover, if $\mathbf{s} \neq \boldsymbol{s}_{\mathcal{M}}(\mathbf{p})$, then
    \begin{equation}\label{eq:strict_conv}
        \mathbf{p}(\boldsymbol{s}_{\mathcal{M}}(\mathbf{p}) - \mathbf{s}) > 0.
    \end{equation}
    
    \begin{lem}\label{lem:cont_pos_ord}
        $\boldsymbol{s}_{\mathcal{M}}(\alpha\mathbf{p}) = \boldsymbol{s}_{\mathcal{M}}(\mathbf{p})$ for all $\alpha \in \mathbb{R}^+$, $\mathbf{p} \in \mathcal{S}^*$, $\mathbf{p} \neq \boldsymbol{0}$.
    \end{lem}
    \begin{pf}
        It follows directly from~\eqref{eq:cont_def}.
    \end{pf}
    Using Lemma~\ref{lem:cont_pos_ord}, we can rewrite~\eqref{eq:dM} by the following way:
    \begin{equation}\label{eq:dM_norm}
        \partial\mathcal{M} = \{\boldsymbol{s}_{\mathcal{M}}(\mathbf{p}):\: \mathbf{p} \in \mathcal{S}^*,\: \lVert\mathbf{p}\rVert = 1\}.
    \end{equation}

    \begin{lem}\label{lem:cont_cont}
        $\boldsymbol{s}_{\mathcal{M}}$ is a continuous function on $\mathcal{S}^* \setminus \{\boldsymbol{0}\}$.
    \end{lem}
    \begin{pf}
        Suppose a contrary, that there exists a sequence $\{\mathbf{p}_n\}_{n = 1}^{\infty}$, such that $\lim_{n \to \infty}\mathbf{p}_n = \mathbf{p}$, and
        \begin{equation}\label{eq:help_eps}
            \lVert\boldsymbol{s}_{\mathcal{M}}(\mathbf{p}) - \boldsymbol{s}_{\mathcal{M}}(\mathbf{p}_n)\rVert \geq \varepsilon
        \end{equation}
        for some $\varepsilon \in \mathbb{R}^+$. Consider a sequence $\mathbf{s}_n = \boldsymbol{s}_{\mathcal{M}}(\mathbf{p}_n) \in \mathcal{M}$. Since $\mathcal{M}$ is a compact set, there exists sub-sequence $\{\mathbf{s}_{n_k}\}_{k = 1}^{\infty}$ such that
        \begin{equation*}
            \lim_{k \to \infty}\mathbf{s}_{n_k} = \tilde{\mathbf{s}} \in \mathcal{M}.
        \end{equation*}
        Note, that $\tilde{\mathbf{s}} \neq \boldsymbol{s}_{\mathcal{M}}(\mathbf{p})$, because
        \begin{equation*}
            \lVert\boldsymbol{s}_{\mathcal{M}}(\mathbf{p}) - \tilde{\mathbf{s}}\rVert = \lim_{k \to \infty}\lVert\boldsymbol{s}_{\mathcal{M}}(\mathbf{p}) - \boldsymbol{s}_{\mathcal{M}}(\mathbf{p}_{n_k})\rVert \geq \varepsilon > 0.
        \end{equation*}
        Using~\eqref{eq:conv} for $\mathbf{p} = \mathbf{p}_{n_k}$, $\mathbf{s} = \boldsymbol{s}_{\mathcal{M}}(\mathbf{p})$, we obtain
        \begin{equation*}
            \mathbf{p}_{n_k}(\boldsymbol{s}_{\mathcal{M}}(\mathbf{p}_{n_k}) - \boldsymbol{s}_{\mathcal{M}}(\mathbf{p})) \geq 0.
        \end{equation*}
        Thus,
        \begin{equation*}
            \mathbf{p}(\tilde{\mathbf{s}} - \boldsymbol{s}_{\mathcal{M}}(\mathbf{p})) = \lim_{k \to \infty} \mathbf{p}(\boldsymbol{s}_{\mathcal{M}}(\mathbf{p}_{n_k}) - \boldsymbol{s}_{\mathcal{M}}(\mathbf{p})) \geq 0.
        \end{equation*}
        On the other hand, using~\eqref{eq:strict_conv}, we obtain
        \begin{equation*}
            \mathbf{p}(\boldsymbol{s}_{\mathcal{M}}(\mathbf{p}) - \tilde{\mathbf{s}}) > 0.
        \end{equation*}
        The contradiction completes the proof.
    \end{pf}

    For the next property of contact functions we use the notions of the directional derivative and gradient. If $f: \mathcal{S}^* \to \mathbb{R}$, then the gradient $\frac{\partial f}{\partial\mathbf{p}}(\mathbf{p}) \in \mathcal{S}$. Let $\mathbf{q} \in \mathcal{S}^*$ be a direction, then the directional derivative is given by
    \begin{equation*}
        \nabla_{\mathbf{p}, \mathbf{q}}f(\mathbf{p}) \overset{\mathrm{def}}{=} \lim_{\gamma \to 0} \frac{f(\mathbf{p} + \gamma\mathbf{q}) - f(\mathbf{p})}{\gamma} = \mathbf{q}\frac{\partial f}{\partial\mathbf{p}}(\mathbf{p}).
    \end{equation*}
    \begin{lem}\label{lem:cont_p_diff}
        If $\mathbf{p} \in \mathcal{S}^*$, $\mathbf{p} \neq \boldsymbol{0}$, then
        \begin{equation*}
            \frac{\partial}{\partial\mathbf{p}}(\mathbf{p}\boldsymbol{s}_{\mathcal{M}}(\mathbf{p})) = \boldsymbol{s}_{\mathcal{M}}(\mathbf{p}).
        \end{equation*}
    \end{lem}
    \begin{pf}
        At first we prove that
        \begin{equation}\label{eq:tmp_lim}
            \lim_{\gamma \to 0} \mathbf{p}\frac{\boldsymbol{s}_{\mathcal{M}}(\mathbf{p} + \gamma\mathbf{q}) - \boldsymbol{s}_{\mathcal{M}}(\mathbf{p})}{\gamma} = 0.
        \end{equation}
        Using the Schwarz inequality and~\eqref{eq:conv}, we obtain
        \begin{multline*}
            \lVert\gamma\mathbf{q}\rVert\lVert\boldsymbol{s}_{\mathcal{M}}(\mathbf{p} + \gamma\mathbf{q}) - \boldsymbol{s}_{\mathcal{M}}(\mathbf{p})\rVert\\
            \geq \gamma\mathbf{q}(\boldsymbol{s}_{\mathcal{M}}(\mathbf{p} + \gamma\mathbf{q}) - \boldsymbol{s}_{\mathcal{M}}(\mathbf{p}))\\
            = (\mathbf{p} + \gamma\mathbf{q})(\boldsymbol{s}_{\mathcal{M}}(\mathbf{p} + \gamma\mathbf{q}) - \boldsymbol{s}_{\mathcal{M}}(\mathbf{p}))\\
            + \mathbf{p}(\boldsymbol{s}_{\mathcal{M}}(\mathbf{p}) - \boldsymbol{s}_{\mathcal{M}}(\mathbf{p} + \gamma\mathbf{q}))\\
            \geq \mathbf{p}(\boldsymbol{s}_{\mathcal{M}}(\mathbf{p}) - \boldsymbol{s}_{\mathcal{M}}(\mathbf{p} + \gamma\mathbf{q})) \geq 0.
        \end{multline*}
        Dividing by $\gamma$ and using Lemma~\ref{lem:cont_cont}, we obtain~\eqref{eq:tmp_lim}. Thus,
        \begin{multline*}
            \nabla_{\mathbf{p}, \mathbf{q}}(\mathbf{p}\boldsymbol{s}_{\mathcal{M}}(\mathbf{p})) = \lim_{\gamma \to 0} \frac{(\mathbf{p} + \gamma\mathbf{q})\boldsymbol{s}_{\mathcal{M}}(\mathbf{p} + \gamma\mathbf{q}) - \mathbf{p}\boldsymbol{s}_{\mathcal{M}}(\mathbf{p})}{\gamma}\\
            = \lim_{\gamma \to 0} \mathbf{q}\boldsymbol{s}_{\mathcal{M}}(\mathbf{p} + \gamma\mathbf{q}) + \lim_{\gamma \to 0} \mathbf{p}\frac{\boldsymbol{s}_{\mathcal{M}}(\mathbf{p} + \gamma\mathbf{q}) - \boldsymbol{s}_{\mathcal{M}}(\mathbf{p})}{\gamma}\\
            = \mathbf{q}\boldsymbol{s}_{\mathcal{M}}(\mathbf{p}).
        \end{multline*}
        The proof is completed by linking the directional derivative with the gradient.
    \end{pf}

    The following property will also be extremely useful to us.
    \begin{lem}\label{lem:MNDT}
        Let $\mathcal{M}_1$, $\mathcal{M}_2$ be strictly convex compact sets that intersect at most one point: $\mathrm{card}(\mathcal{M}_1 \cap \mathcal{M}_2) \leq 1$. Then
        \begin{multline*}
            \min_{\mathbf{s}_1 \in \mathcal{M}_1,\:\mathbf{s}_2 \in \mathcal{M}_2}\lVert\mathbf{s}_1 - \mathbf{s}_2\rVert \\
            = \max_{\mathbf{p} \in \mathcal{S}^*,\: \lVert\mathbf{p}\rVert = 1}\mathbf{p}(\boldsymbol{s}_{\mathcal{M}_2}(-\mathbf{p}) - \boldsymbol{s}_{\mathcal{M}_1}(\mathbf{p})).
        \end{multline*}
    \end{lem}
    This statement is known as Minimum Norm Duality Theorem~\citep{Dax2006-iq}.
    
    \section{Problem formulation}\label{sec:prob_st}

    \subsection{General problem}
    
    The dynamics of the system is described by a linear differential system
    \begin{equation}\label{eq:dyn}
        \dot{\boldsymbol{s}} = \boldsymbol{A}\boldsymbol{s} + \boldsymbol{u}.
    \end{equation}
    Here, the state-vector $\boldsymbol{s}(t; \boldsymbol{u}) \in \mathcal{S}$ at time $t \in \mathbb{R}^+_0$ corresponds to an absolutely continuous solution $\boldsymbol{s}(\cdot; \boldsymbol{u})$ of~\eqref{eq:dyn} for a fixed control input $\boldsymbol{u}$ with the initial conditions $\boldsymbol{s}(0; \boldsymbol{u}) = \mathbf{s}_0 \in \mathcal{S}$. We suppose that the control input is restricted $\boldsymbol{u}(t) \in \mathcal{U} \subset \mathcal{S}$, where $\mathcal{U}$ is a compact set. Each measurable function with range of values $\mathcal{U}$ is an admissible control input. The set of all admissible control inputs is denoted by $\mathcal{A}$. Thus, any $\boldsymbol{u} \in \mathcal{A}$ is such that for all $t \in \mathbb{R}^+_0$: $\boldsymbol{u}(t) \in \mathcal{U}$.
    
    $\boldsymbol{A}: \mathbb{R} \to \mathbb{R}^{n\times n}$ is a matrix-valued function that is measurable on $\mathbb{R}$. We denote the fundamental solution matrix by $\boldsymbol{\Phi}(t) \in \mathbb{R}^{n \times n}$, i.e. $\dot{\boldsymbol{\Phi}} = \boldsymbol{A}\boldsymbol{\Phi}$, $\boldsymbol{\Phi}(0) = \mathbf{I}$\footnote{$\mathbf{I}$ is the identity matrix in $\mathbb{R}^{n \times n}$.}. Thus,
    \begin{equation}\label{eq:sol_through_fund}
        \boldsymbol{s}(t; \boldsymbol{u}) = \boldsymbol{\Phi}(t)\mathbf{s}_0 + \boldsymbol{\Phi}(t)\int_0^t\boldsymbol{\Phi}(\tau)^{-1}\boldsymbol{u}(\tau)\mathrm{d}\tau.
    \end{equation}
    Note, that if $\boldsymbol{A}(t) = \mathbf{A} = \mathrm{const}$, then $\boldsymbol{\Phi}(t) = e^{t\mathbf{A}}$. The fundamental solution matrix is continuous and has continuous derivatives. This matrix is non-singular for all $t \in \mathbb{R}^+_0$.

    The target set is given by time-varying set $\mathcal{G}(t) \subset \mathcal{S}$, i.e. $\mathcal{G}:\mathbb{R} \to 2^{\mathcal{S}}$ is a set-valued mapping. We define the general problem by finding a control input $\boldsymbol{u} \in \mathcal{A}$ that sets the minimum value of the following functional:
    \begin{equation*}
        J[\boldsymbol{u}] \overset{\mathrm{def}}{=} \inf \left\{t \in \mathbb{R}^+_0:\: \boldsymbol{s}(t; \boldsymbol{u}) \in \mathcal{G}(t)\right\}.
    \end{equation*}
    The minimum-time required to reach the target set is defined by
    \begin{equation*}
        T^* \overset{\mathrm{def}}{=} \inf_{\boldsymbol{u} \in \mathcal{A}} J[\boldsymbol{u}].
    \end{equation*} 
    If the target set cannot be visited we formally suppose that $T^* = +\infty$. Our goal is to describe algorithms for computing $T^*$ under some restrictive assumptions that are described further.



    \subsection{Assumptions}

    The general problem is quite difficult for computational algorithms. Now we formulate two groups of restrictions, using which we can construct convergent algorithms for solving the problem.

    Let $\mathcal{R}(t)$ be a reachable set at time $t \in \mathbb{R}^+_0$, i.e. it is a set of ending points of trajectories that can be raised by all admissible control inputs $\boldsymbol{u} \in \mathcal{A}$:
    \begin{equation*}
        \mathcal{R}(t) \overset{\mathrm{def}}{=} \left\{\boldsymbol{s}(t; \boldsymbol{u}):\: \boldsymbol{u} \in \mathcal{A}\right\}.
    \end{equation*}
    The first group of restrictions relates to the plant:
    \begin{enumerate}
        \item[$\mathrm{RP1}$] (normality) For all $\boldsymbol{u}_1 \in \mathcal{A}$, $\boldsymbol{u}_2 \in \mathcal{A}$ if $\boldsymbol{s}(t; \boldsymbol{u}_1) = \boldsymbol{s}(t; \boldsymbol{u}_2) \in \partial\mathcal{R}(t)$ then $\boldsymbol{u}_1(\tau) \overset{\mathrm{a.e.}}{=} \boldsymbol{u}_2(\tau)$, $\tau \in [0, t]$.
        \item[$\mathrm{RP2}$] (speed limit) There exists $v_{\mathcal{R}} \in \mathbb{R}^+$ such that for any $\boldsymbol{u} \in \mathcal{A}$ and $t \in \mathbb{R}^+_0$: $\lVert\dot{\boldsymbol{s}}(t; \boldsymbol{u})\rVert = \lVert\boldsymbol{A}(t)\boldsymbol{s}(t; \boldsymbol{u}) + \boldsymbol{u}(t)\rVert \leq v_{\mathcal{R}}$.
    \end{enumerate}
    
    The second group of restrictions is devoted to the target set:
    \begin{enumerate}
        \item[$\mathrm{RT1}$] (simple geometry) $\mathcal{G}(t)$ is a strictly convex compact set for any $t \in \mathbb{R}$.
        \item[$\mathrm{RT2}$] (continuity) $\mathcal{G}$ is a continuous mapping and the contact function $\boldsymbol{s}_{\mathcal{G}(\cdot)}(\mathbf{p})$ is continuous for all $\mathbf{p} \in \mathcal{S}^*$, $\mathbf{p} \neq \boldsymbol{0}$ ($\mathrm{RT1}$ holds).
        \item[$\mathrm{RT3}$] (differentiability) $\boldsymbol{s}_{\mathcal{G}(\cdot)}(\mathbf{p})$ is differentiable for all $\mathbf{p} \in \mathcal{S}^*$, $\mathbf{p} \neq \boldsymbol{0}$ and $\boldsymbol{v}_{\mathcal{G}(t)}(\mathbf{p}) = \frac{\partial}{\partial t} \boldsymbol{s}_{\mathcal{G}(t)}(\mathbf{p})$ is known ($\mathrm{RT2}$ holds).
        \item[$\mathrm{RT4}$] (speed limit) $\lVert\boldsymbol{v}_{\mathcal{G}(t)}(\mathbf{p})\rVert \leq v_{\mathcal{G}} \in \mathbb{R}^+_0$ for all $\mathbf{p} \in \mathcal{S}^*$, $\mathbf{p} \neq \boldsymbol{0}$ ($\mathrm{RT3}$ holds).
    \end{enumerate}

    \section{Geometric view}\label{sec:constr}

    In this section, we introduce a series of fundamental concepts that provide geometric representations of actions and are utilized in constructing algorithms for solving the general problem.
    
    \subsection{Maximum principle}
    
    The adjoint system is described by a linear differential system
    \begin{equation}\label{eq:adj}
        \dot{\boldsymbol{p}} = -\boldsymbol{p}\boldsymbol{A}, \quad \boldsymbol{p}(T; T, \mathbf{p}) = \mathbf{p} \in \mathcal{S}^*, \quad \mathbf{p} \neq \boldsymbol{0}.
    \end{equation}
    We associate $T \in \mathbb{R}^+_0$ with some terminal time moment when the solution $\boldsymbol{p}(\cdot; T, \mathbf{p})$ of~\eqref{eq:adj} is equal to $\mathbf{p} \in \mathcal{S}^*$. Using the fundamental solution matrix, we obtain $\boldsymbol{p}(t; T, \mathbf{p}) = \mathbf{p}\boldsymbol{\Phi}(T)\boldsymbol{\Phi}(t)^{-1}$.

    According to the maximum principle for linear systems, the extremal control inputs must satisfy the maximum condition:
    \begin{equation}\label{eq:max_princ}
        \boldsymbol{p}(t; T, \mathbf{p})\boldsymbol{u}(t) \overset{\mathrm{a.e.}}{=} \max_{\mathbf{u} \in \mathcal{U}} \boldsymbol{p}(t; T, \mathbf{p})\mathbf{u}.
    \end{equation}
    If $\mathrm{RP1}$ holds, then Theorem 3 of~\citet[pp.~76--77]{Lee1967-vf} guarantees that~\eqref{eq:max_princ} determines the extremal control input $\boldsymbol{u}_E(\boldsymbol{p}(\cdot; T, \mathbf{p}))$ for almost all $t \in [0, T]$. Moreover, $\mathcal{R}(t)$ is a strictly convex compact set.
    
    We suppose that $\boldsymbol{u}_E$ can be easily expressed. For example, if $\mathcal{U}$ is a unit ball, i.e. $\mathcal{U} = \{\mathbf{u} \in \mathcal{S}:\: \lVert\mathbf{u}\rVert \leq 1\}$, then $\boldsymbol{u}_E(\mathbf{p}) = \mathbf{p}^\top/\lVert\mathbf{p}\rVert$.

    \subsection{Boundary of the reachable set}

    Using~\eqref{eq:sol_through_fund} and $\boldsymbol{u}_E$, we can describe extremal trajectories that lead to the boundary of the reachable set:
    \begin{multline*}
        \boldsymbol{s}_E(t; T, \mathbf{p}) \overset{\mathrm{def}}{=} \boldsymbol{s}(t; \boldsymbol{u}_E(\boldsymbol{p}(\cdot; T, \mathbf{p}))) \\
        = \boldsymbol{\Phi}(t)\mathbf{s}_0 + \boldsymbol{\Phi}(t)\int_0^t \boldsymbol{\Phi}(\tau)^{-1}\boldsymbol{u}_E(\boldsymbol{p}(\tau; T, \mathbf{p}))\mathrm{d}\tau.
    \end{multline*}
    These trajectories are parameterized by $\mathbf{p} \in \mathcal{S}^*$ ($\mathbf{p} \neq \boldsymbol{0}$).

\begin{figure}
        \begin{center}
            \includegraphics[width=0.4\textwidth]{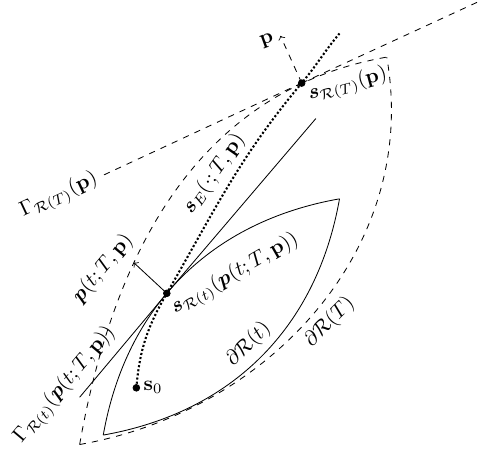}
            \caption{Connection of extremal trajectory $\boldsymbol{s}_E(\cdot; T, \mathbf{p})$ and contact function $\boldsymbol{s}_{\mathcal{R}(\cdot)}(\boldsymbol{p}(\cdot; T, \mathbf{p}))$. Here, $t < T$.}
            \label{fig:traj_ex}
        \end{center}
    \end{figure}

    The point $\boldsymbol{s}_E(t; T, \mathbf{p})$ is a terminal point of extremal trajectory, i.e. this point lies on $\partial\mathcal{R}(t)$. Moreover, $\boldsymbol{p}(t; T, \mathbf{p})$ is an outward normal to $\mathcal{R}(t)$ at point $\boldsymbol{s}_E(t; T, \mathbf{p})$~\citep[Corollary 1, p.~75]{Lee1967-vf}. Thus, $\boldsymbol{s}_E(t; T, \mathbf{p})$ is a contact point that corresponds to the support vector $\boldsymbol{p}(t; T, \mathbf{p})$, i.e. $\boldsymbol{s}_{\mathcal{R}(t)}(\boldsymbol{p}(t; T, \mathbf{p})) = \boldsymbol{s}_E(t; T, \mathbf{p})$\footnote{We use the notion $\boldsymbol{s}_{\mathcal{M}}$ for $\mathcal{M} = \mathcal{R}(t)$ from Section~\ref{sec:prel}} (see~Fig.~\ref{fig:traj_ex}).

    Using~$\boldsymbol{p}(t; t, \mathbf{p}) = \mathbf{p}$ and~\eqref{eq:dM_norm}, we obtain the parametric description of the boundary of the reachable set:
    \begin{equation*}
        \partial\mathcal{R}(t) = \{\boldsymbol{s}_E(t; t, \mathbf{p}):\: \mathbf{p} \in \mathcal{S}^*, \: \lVert\mathbf{p}\rVert = 1\}.
    \end{equation*}

    \begin{lem}\label{lem:sR_cont_t}
         \textup{[\citet{Neustadt1960-ki}]}. If $\mathrm{RP1}$ holds than $\boldsymbol{s}_{\mathcal{R}(\cdot)}(\mathbf{p})$ is a continuous function for any $\mathbf{p} \in \mathcal{S}^*$, $\mathbf{p} \neq \boldsymbol{0}$.
    \end{lem}
    \begin{pf}
        Lemma~\ref{lem:cont_cont} and $\boldsymbol{p}(t; t, \mathbf{p}) = \mathbf{p}$ give
        \begin{equation*}
            \lim_{\theta \to 0} \boldsymbol{s}_{\mathcal{R}(t)}(\boldsymbol{p}(t; t + \theta, \mathbf{p})) = \boldsymbol{s}_{\mathcal{R}(t)}(\mathbf{p}).
        \end{equation*}
        Thus,
        \begin{equation*}
            \lim_{\theta \to 0}\lVert\boldsymbol{s}_{\mathcal{R}(t)}(\boldsymbol{p}(t; t + \theta, \mathbf{p})) - \boldsymbol{s}_{\mathcal{R}(t)}(\mathbf{p})\rVert = 0.
        \end{equation*}
        Using the definition of $\boldsymbol{s}_E$, we deduce
        \begin{multline*}
            \boldsymbol{s}_E(t + \theta; t + \theta, \mathbf{p}) - \boldsymbol{s}_E(t; t + \theta, \mathbf{p}) = (\boldsymbol{\Phi}(t + \theta) - \boldsymbol{\Phi}(t))\mathbf{s}_0\\
            + (\boldsymbol{\Phi}(t + \theta) - \boldsymbol{\Phi}(t))\int_0^t \boldsymbol{\Phi}(\tau)^{-1}\boldsymbol{u}_E(\boldsymbol{p}(\tau; t + \theta, \mathbf{p}))\mathrm{d}\tau\\
            + \boldsymbol{\Phi}(t + \theta)\int_t^{t + \theta} \boldsymbol{\Phi}(\tau)^{-1}\boldsymbol{u}_E(\boldsymbol{p}(\tau; t + \theta, \mathbf{p}))\mathrm{d}\tau.
        \end{multline*}
        The continuity of $\boldsymbol{\Phi}$ and boundedness of $\boldsymbol{u}_E$ give
        \begin{equation*}
            \lim_{\theta \to 0}\lVert\boldsymbol{s}_E(t + \theta; t + \theta, \mathbf{p}) - \boldsymbol{s}_E(t; t + \theta, \mathbf{p})\rVert = 0.
        \end{equation*}
        Finally, using $\boldsymbol{s}_E(t + \theta; t + \theta, \mathbf{p}) = \boldsymbol{s}_{\mathcal{R}(t + \theta)}(\mathbf{p})$ and $\boldsymbol{s}_E(t; t + \theta, \mathbf{p}) = \boldsymbol{s}_{\mathcal{R}(t)}(\boldsymbol{p}(t; t + \theta, \mathbf{p}))$, we deduce
        \begin{multline*}
            \lim_{\theta \to 0}\lVert\boldsymbol{s}_{\mathcal{R}(t + \theta)}(\mathbf{p}) - \boldsymbol{s}_{\mathcal{R}(t)}(\mathbf{p})\rVert\\
            \leq \lim_{\theta \to 0}\lVert\boldsymbol{s}_{\mathcal{R}(t + \theta)}(\mathbf{p}) - \boldsymbol{s}_{\mathcal{R}(t)}(\boldsymbol{p}(t; t + \theta, \mathbf{p}))\rVert\\
            + \lim_{\theta \to 0}\lVert\boldsymbol{s}_{\mathcal{R}(t)}(\boldsymbol{p}(t; t + \theta, \mathbf{p})) - \boldsymbol{s}_{\mathcal{R}(t)}(\mathbf{p})\rVert = 0.
        \end{multline*}
    \end{pf}
    
    \begin{lem}\label{lem:psR_dt}
        Let $\mathrm{RP1}$, $\mathbf{p} \in \mathcal{S}^*$, $\mathbf{p} \neq 0$. Then
        \begin{equation*}
            \frac{\mathrm{d}}{\mathrm{d}t}(\mathbf{p}\boldsymbol{s}_{\mathcal{R}(t)}(\mathbf{p})) = \mathbf{p}\boldsymbol{A}(t)\boldsymbol{s}_{\mathcal{R}(t)}(\mathbf{p}) + \mathbf{p}\boldsymbol{u}_E(\mathbf{p}).
        \end{equation*}
    \end{lem}
    \begin{pf}
        Note, that
        \begin{equation*}
            \lim_{\theta \to 0}\mathbf{p}\frac{\boldsymbol{s}_{\mathcal{R}(t + \theta)}(\mathbf{p}) - \boldsymbol{s}_{\mathcal{R}(t + \theta)}(\boldsymbol{p}(t + \theta; t, \mathbf{p}))}{\theta} = 0.
        \end{equation*}
        can be evaluated in the same way as~\eqref{eq:tmp_lim}. Using $\boldsymbol{s}_{\mathcal{R}(t + \theta)}(\boldsymbol{p}(t + \theta; t, \mathbf{p})) = \boldsymbol{s}_E(t + \theta; t, \mathbf{p})$, we have
        \begin{multline*}
            \frac{\mathrm{d}}{\mathrm{d}t}(\mathbf{p}\boldsymbol{s}_{\mathcal{R}(t)}(\mathbf{p})) = \lim_{\theta \to 0}\mathbf{p}\frac{\boldsymbol{s}_{\mathcal{R}(t + \theta)}(\mathbf{p}) - \boldsymbol{s}_{\mathcal{R}(t)}(\mathbf{p})}{\theta}\\
            = \lim_{\theta \to 0}\mathbf{p}\frac{\boldsymbol{s}_{\mathcal{R}(t + \theta)}(\boldsymbol{p}(t + \theta; t, \mathbf{p})) - \boldsymbol{s}_{\mathcal{R}(t)}(\mathbf{p})}{\theta}\\
            = \lim_{\theta \to 0}\mathbf{p}\frac{\boldsymbol{s}_E(t + \theta; t, \mathbf{p}) - \boldsymbol{s}_E(t; t, \mathbf{p})}{\theta}\\
            = \mathbf{p}\boldsymbol{A}(t)\boldsymbol{s}_{\mathcal{R}(t)}(\mathbf{p}) + \mathbf{p}\boldsymbol{u}_E(\mathbf{p}).
        \end{multline*}
    \end{pf}

    \subsection{Distance estimates}
    
    Now we describe some properties of function that corresponds to the distance from the target set $\mathcal{G}(t)$ to the reachable set $\mathcal{R}(t)$:
    \begin{equation*}
        \rho(t) \overset{\mathrm{def}}{=} \min_{\mathbf{s} \in \mathcal{R}(t),\:\tilde{\mathbf{s}} \in \mathcal{G}(t)}\lVert\tilde{\mathbf{s}} - \mathbf{s}\rVert.
    \end{equation*}
    The distance function helps to express the minimum-time required to reach the target set:
    \begin{equation}\label{eq:T_def}
        T^* = \min\{t \in \mathbb{R}^+_0: \rho(t) = 0\}.
    \end{equation}
    Here we suppose that $\min\varnothing = +\infty$. If the target set can be visited, then there exists $t \in \mathbb{R}^+_0$ such that the distance vanishes and $T^* < +\infty$. Otherwise, we technically suppose that $T^* = +\infty$. 

\begin{figure}
        \begin{center}
            \includegraphics[width=0.3\textwidth]{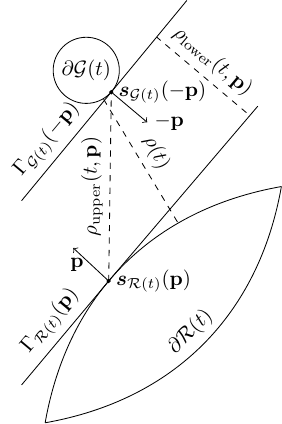}
            \caption{Distance estimates $\rho_{\mathrm{lower}}(t, \mathbf{p}) \leq \rho(t) \leq \rho_{\mathrm{upper}}(t, \mathbf{p})$}
            \label{fig:estim}
        \end{center}
    \end{figure}

    Using Lemma~\ref{lem:MNDT}, we obtain
    \begin{equation}\label{eq:MNDT}
        \rho(t) = \max_{\mathbf{p} \in \mathcal{S}^*, \: \lVert\mathbf{p}\rVert = 1}\mathbf{p}(\boldsymbol{s}_{\mathcal{G}(t)}(-\mathbf{p}) - \boldsymbol{s}_{\mathcal{R}(t)}(\mathbf{p}))
    \end{equation}
    for $\mathrm{card}(\mathcal{G}(t) \cap \mathcal{R}(t)) \leq 1$. Lemma~\ref{lem:cont_pos_ord} gives
    \begin{multline}\label{eq:lower}
        \rho(t) = \max_{\mathbf{p} \in \mathcal{S}^*, \: \mathbf{p} \neq \boldsymbol{0}}\frac{\mathbf{p}}{\lVert\mathbf{p}\rVert}(\boldsymbol{s}_{\mathcal{G}(t)}(-\mathbf{p}) - \boldsymbol{s}_{\mathcal{R}(t)}(\mathbf{p}))\\
        \geq \frac{\mathbf{p}}{\lVert\mathbf{p}\rVert}(\boldsymbol{s}_{\mathcal{G}(t)}(-\mathbf{p}) - \boldsymbol{s}_{\mathcal{R}(t)}(\mathbf{p}))
    \end{multline}
    for any $\mathbf{p} \in \mathcal{S}^*$, $\mathbf{p} \neq \boldsymbol{0}$. Therefore, the value 
    \begin{equation*}
        \rho_{\mathrm{lower}}(t, \mathbf{p}) \overset{\mathrm{def}}{=} \frac{\mathbf{p}}{\lVert\mathbf{p}\rVert}(\boldsymbol{s}_{\mathcal{G}(t)}(-\mathbf{p}) - \boldsymbol{s}_{\mathcal{R}(t)}(\mathbf{p}))
    \end{equation*}
    is a lower estimate of $\rho(t)$ (see~Fig.~\ref{fig:estim}). Note, that $|\rho_{\mathrm{lower}}(t, \mathbf{p})|$ corresponds to the distance between $\Gamma_{\mathcal{R}(t)}(\mathbf{p})$ and $\Gamma_{\mathcal{G}(t)}(-\mathbf{p})$. These hyperplanes are tangent hyperplanes for $\mathcal{R}(t)$ and $\mathcal{G}(t)$ at points $\boldsymbol{s}_{\mathcal{R}(t)}(\mathbf{p})$ and $\boldsymbol{s}_{\mathcal{G}(t)}(-\mathbf{p})$.
    
    On the other hand, we have
    \begin{equation}\label{eq:upper}
        \rho(t) = \min_{\mathbf{s} \in \mathcal{R}(t),\:\tilde{\mathbf{s}} \in \mathcal{G}(t)}\lVert\tilde{\mathbf{s}} - \mathbf{s}\rVert \leq \lVert\boldsymbol{s}_{\mathcal{G}(t)}(-\mathbf{p}) - \boldsymbol{s}_{\mathcal{R}(t)}(\mathbf{p})\rVert
    \end{equation}
    for any $\mathbf{p} \in \mathcal{S}^*$, $\mathbf{p} \neq \boldsymbol{0}$. Hence, the value
    \begin{equation*}
        \rho_{\mathrm{upper}}(t, \mathbf{p}) \overset{\mathrm{def}}{=} \lVert\boldsymbol{s}_{\mathcal{G}(t)}(-\mathbf{p}) - \boldsymbol{s}_{\mathcal{R}(t)}(\mathbf{p})\rVert
    \end{equation*}
    is an upper estimate of $\rho(t)$ (see~Fig.~\ref{fig:estim}).

    \begin{lem}\label{lem:upper_lower}
        Let $\mathrm{RP1}$, $\mathrm{RT1}$, $\mathbf{p}, \tilde{\mathbf{p}} \in \mathcal{S}^*$, $\mathbf{p} \neq \boldsymbol{0}$, $\tilde{\mathbf{p}} \neq \boldsymbol{0}$. Then
        \begin{equation*}
            \rho_{\mathrm{lower}}(t, \mathbf{p}) \leq \rho_{\mathrm{upper}}(t, \tilde{\mathbf{p}}).
        \end{equation*}
    \end{lem}
    \begin{pf}
        It follows directly from~\eqref{eq:lower} and~\eqref{eq:upper} if $\mathrm{card}(\mathcal{G}(t) \cap \mathcal{R}(t)) \leq 1$. If $\mathrm{card}(\mathcal{G}(t) \cap \mathcal{R}(t)) > 1$, then $\mathcal{G}(t)$, $\mathcal{R}(t)$ are not separable. Thus, $\rho_{\mathrm{lower}}(t, \mathbf{p}) < 0 \leq \rho_{\mathrm{upper}}(t, \tilde{\mathbf{p}})$.
    \end{pf}

    \begin{lem}\label{lem:lower_diff_p}
        \textup{[\citet{Gilbert1966-sc}]}. Let $\mathrm{RP1}$, $\mathrm{RT1}$, $\mathbf{p} \in \mathcal{S}^*$, $\mathbf{p} \neq \boldsymbol{0}$. Then
        \begin{equation*}
            \frac{\partial\rho_{\mathrm{lower}}}{\partial\mathbf{p}}(t, \mathbf{p}) = \frac{\boldsymbol{s}_{\mathcal{G}(t)}(-\mathbf{p}) - \boldsymbol{s}_{\mathcal{R}(t)}(\mathbf{p})}{\lVert\mathbf{p}\rVert} - \frac{\mathbf{p}^\top}{\lVert\mathbf{p}\rVert^2}\rho_{\mathrm{lower}}(t, \mathbf{p}).
        \end{equation*}
    \end{lem}
    \begin{pf}
        Using Lemma~\ref{lem:cont_p_diff}, we obtain
        \begin{equation}\label{eq:xi_will}
            \frac{\partial}{\partial\mathbf{p}}(\lVert\mathbf{p}\rVert\rho_{\mathrm{lower}}(t, \mathbf{p})) = \boldsymbol{s}_{\mathcal{G}(t)}(-\mathbf{p}) - \boldsymbol{s}_{\mathcal{R}(t)}(\mathbf{p}).
        \end{equation}
        Taking into account $\frac{\partial}{\partial\mathbf{p}}\frac{1}{\lVert\mathbf{p}\rVert} = - \frac{\mathbf{p}^\top}{\lVert\mathbf{p}\rVert^3}$, we have
        \begin{multline*}
            \frac{\partial\rho_{\mathrm{lower}}}{\partial\mathbf{p}}(t, \mathbf{p}) = \frac{\partial}{\partial\mathbf{p}}\left(\lVert\mathbf{p}\rVert\rho_{\mathrm{lower}}(t, \mathbf{p})\cdot\frac{1}{\lVert\mathbf{p}\rVert}\right)\\
            = (\boldsymbol{s}_{\mathcal{G}(t)}(-\mathbf{p}) - \boldsymbol{s}_{\mathcal{R}(t)}(\mathbf{p}))\cdot\frac{1}{\lVert\mathbf{p}\rVert} - \lVert\mathbf{p}\rVert\rho_{\mathrm{lower}}(t, \mathbf{p})\cdot\frac{\mathbf{p}^\top}{\lVert\mathbf{p}\rVert^3}\\
            = \frac{\boldsymbol{s}_{\mathcal{G}(t)}(-\mathbf{p}) - \boldsymbol{s}_{\mathcal{R}(t)}(\mathbf{p})}{\lVert\mathbf{p}\rVert} - \frac{\mathbf{p}^\top}{\lVert\mathbf{p}\rVert^2}\rho_{\mathrm{lower}}(t, \mathbf{p}).
        \end{multline*}
    \end{pf}
    
    \begin{lem}\label{lem:lower_diff_t}
        Let $\mathrm{RP1}$, $\mathrm{RT3}$, $\mathbf{p} \in \mathcal{S}^*$, $\mathbf{p} \neq \boldsymbol{0}$. Then
        \begin{multline*}
            \frac{\partial\rho_{\mathrm{lower}}}{\partial t}(t, \mathbf{p})\\
            = \frac{\mathbf{p}}{\lVert\mathbf{p}\rVert}(\boldsymbol{v}_{\mathcal{G}(t)}(-\mathbf{p}) - \boldsymbol{A}(t)\boldsymbol{s}_{\mathcal{R}(t)}(\mathbf{p}) - \boldsymbol{u}_E(\mathbf{p})).
        \end{multline*}
    \end{lem}
    \begin{pf}
        It follows from Lemma~\ref{lem:psR_dt} and RT3.
    \end{pf}

    We denote the difference of upper and lower estimates by the following:
    \begin{equation*}
        \delta(t, \mathbf{p}) \overset{\mathrm{def}}{=} \rho_{\mathrm{upper}}(t, \mathbf{p}) - \rho_{\mathrm{lower}}(t, \mathbf{p}) \geq 0.
    \end{equation*}
    Note that $\delta$ plays a role of the error function. Smaller $\delta(t, \mathbf{p})$ values facilitate greater accuracy in the localization of $\rho(t)$ relative to the upper and lower estimates.
    \begin{lem}\label{lem:delta_cont}
        Let $\mathrm{RP1}$, $\mathrm{RT2}$. Then $\delta$ is a continuous function on $\mathbb{R}^+_0 \times (\mathcal{S}^*\setminus\{\boldsymbol{0}\})$.
    \end{lem}
    \begin{pf}
        It follows from $\mathrm{RT2}$ and Lemmas~\ref{lem:cont_cont},~\ref{lem:sR_cont_t}.
    \end{pf}
    
    If $\delta(t, \mathbf{p}) = 0$, then
    \begin{equation*}
        \rho(t) = \rho_{\mathrm{upper}}(t, \mathbf{p}) = \rho_{\mathrm{lower}}(t, \mathbf{p})  
    \end{equation*}
    and $\boldsymbol{s}_{\mathcal{R}(t)}(\mathbf{p})$ is the nearest point on $\mathcal{R}(t)$ to $\mathcal{G}(t)$, while $\boldsymbol{s}_{\mathcal{G}(t)}(-\mathbf{p})$ is the nearest point on $\mathcal{G}(t)$ to $\mathcal{R}(t)$.

    \begin{lem}\label{lem:min_p_eq}
        Let $\mathrm{RP1}$, $\mathrm{RT1}$, $\mathbf{p} \in \mathcal{S}^*$, $\lVert\mathbf{p}\rVert = 1$. Then, $\delta(t, \mathbf{p}) = 0$ if and only if
        \begin{equation}\label{eq:p_min}
            \mathbf{p} = \left(\frac{\boldsymbol{s}_{\mathcal{G}(t)}(-\mathbf{p}) - \boldsymbol{s}_{\mathcal{R}(t)}(\mathbf{p})}{\lVert\boldsymbol{s}_{\mathcal{G}(t)}(-\mathbf{p}) - \boldsymbol{s}_{\mathcal{R}(t)}(\mathbf{p})\rVert}\right)^\top.
        \end{equation}
    \end{lem}
    \begin{pf}
        The sufficiency is obvious. We prove that $\rho_{\mathrm{lower}}(t, \mathbf{p}) = \rho_{\mathrm{upper}}(t, \mathbf{p})$ implies \eqref{eq:p_min}. Using $\lVert\mathbf{p}\rVert = 1$, we obtain
        \begin{multline*}
            \mathbf{p}(\boldsymbol{s}_{\mathcal{G}(t)}(-\mathbf{p}) - \boldsymbol{s}_{\mathcal{R}(t)}(\mathbf{p})) = \lVert\boldsymbol{s}_{\mathcal{G}(t)}(-\mathbf{p}) - \boldsymbol{s}_{\mathcal{R}(t)}(\mathbf{p})\rVert\\
            = \lVert\mathbf{p}\rVert\lVert\boldsymbol{s}_{\mathcal{G}(t)}(-\mathbf{p}) - \boldsymbol{s}_{\mathcal{R}(t)}(\mathbf{p})\rVert.
        \end{multline*}
        According to Schwarz inequality $\mathbf{p}$ must be co-directed to $\boldsymbol{s}_{\mathcal{G}(t)}(-\mathbf{p}) - \boldsymbol{s}_{\mathcal{R}(t)}(\mathbf{p})$. Thus, \eqref{eq:p_min} holds.
    \end{pf}

    \begin{thm}\label{thm:main_syst}
        Let $\mathrm{RP1}$, $\mathrm{RT1}$, $t \in \mathbb{R}^+_0$, $\mathrm{card}(\mathcal{G}(t) \cap \mathcal{R}(t)) \leq 1$. Then, there exists a unique solution of the system
        \begin{equation*}
            \delta(t, \mathbf{p}) = 0, \quad \lVert\mathbf{p}\rVert = 1, \quad \mathbf{p} \in \mathcal{S}^*.
        \end{equation*}
    \end{thm}
    \begin{pf}
        Using the strict convexity and compactness of $\mathcal{R}(t)$ and $\mathcal{G}(t)$, we conclude that there exists unique $\mathbf{s} \in \partial\mathcal{R}(t)$ and $\tilde{\mathbf{s}} \in \partial\mathcal{G}(t)$ such that $\rho(t) = \lVert\mathbf{s} - \tilde{\mathbf{s}}\rVert$. Using~\eqref{eq:MNDT}, we state that there exists $\mathbf{p} \in \mathcal{S}^*$ ($\lVert\mathbf{p}\rVert = 1$) such that $\rho(t) = \mathbf{p}(\boldsymbol{s}_{\mathcal{G}(t)}(-\mathbf{p}) - \boldsymbol{s}_{\mathcal{R}(t)}(\mathbf{p}))$. Uniqueness of $\mathbf{s}$ and $\tilde{\mathbf{s}}$ guarantees that $\mathbf{s} = \boldsymbol{s}_{\mathcal{R}(t)}(\mathbf{p})$ and $\tilde{\mathbf{s}} = \boldsymbol{s}_{\mathcal{G}(t)}(-\mathbf{p})$. Thus,
        \begin{equation*}
            \lVert\boldsymbol{s}_{\mathcal{G}(t)}(-\mathbf{p}) - \boldsymbol{s}_{\mathcal{R}(t)}(\mathbf{p})\rVert = \mathbf{p}(\boldsymbol{s}_{\mathcal{G}(t)}(-\mathbf{p}) - \boldsymbol{s}_{\mathcal{R}(t)}(\mathbf{p})).
        \end{equation*}
        According to Schwartz inequality, it is fulfilled only if~\eqref{eq:p_min} holds. Lemma~\ref{lem:min_p_eq} completes the proof.
    \end{pf}

    Theorem~\ref{thm:main_syst} declares the existence of a unique solution of
    \begin{equation}\label{eq:p_opt_uniq}
        \delta(T^*, \mathbf{p}) = 0, \quad \lVert\mathbf{p}\rVert = 1, \quad \mathbf{p} \in \mathcal{S}^*
    \end{equation}
    for $T^* < +\infty$. We denote this solution by $\mathbf{p}^*$. Note that $T^*$ and $\mathbf{p}^*$ completely define the optimal control input $\boldsymbol{u}_E(\boldsymbol{p}(\cdot; T^*, \mathbf{p}^*))$. Despite the uniqueness of the solution of~\eqref{eq:p_opt_uniq}, there may exist other vectors $\tilde{\mathbf{p}} \in \mathcal{S}^*$ defining the optimal control input $\boldsymbol{u}_E(\boldsymbol{p}(\cdot; T^*, \tilde{\mathbf{p}}))$. $\mathrm{RP1}$ guarantees that
    \begin{equation*}
        \boldsymbol{u}_E(\boldsymbol{p}(t; T^*, \mathbf{p}^*)) \overset{\mathrm{a.e.}}{=} \boldsymbol{u}_E(\boldsymbol{p}(t; T^*, \tilde{\mathbf{p}})).
    \end{equation*}

    \subsection{Inclination of hyperplanes}

\begin{figure}
        \begin{center}
            \includegraphics[width=0.3\textwidth]{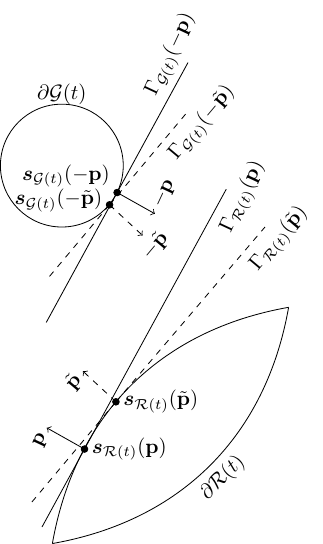}
            \caption{Inclination of hyperplanes $\Gamma_{\mathcal{R}(t)}(\mathbf{p})$, $\Gamma_{\mathcal{G}(t)}(-\mathbf{p})$ by changing $\mathbf{p}$ to $\tilde{\mathbf{p}} = \mathbf{p} + \gamma\boldsymbol{q}(t, \mathbf{p})$. Here, $\gamma$ is such that $\rho_{\mathrm{lower}}(t, \tilde{\mathbf{p}}) > \rho_{\mathrm{lower}}(t, \mathbf{p})$}
            \label{fig:incline}
        \end{center}
    \end{figure}

    Changing the support vector $\mathbf{p}$ for the hyperplane $\Gamma_{\mathcal{M}}(\mathbf{p})$ inclines the hyperplane. If hyperplanes $\Gamma_{\mathcal{R}(t)}(\mathbf{p})$, $\Gamma_{\mathcal{G}(t)}(-\mathbf{p})$ separate sets $\mathcal{R}(t)$, $\mathcal{G}(t)$, but the points $\boldsymbol{s}_{\mathcal{R}(t)}(\mathbf{p})$, $\boldsymbol{s}_{\mathcal{G}(t)}(-\mathbf{p})$ are not closest, i.e. $\rho(t) < \rho_{\mathrm{upper}}(t, \mathbf{p})$, then we can incline these hyperplanes so that the lower bound for the distance increases (see Fig.~\ref{fig:incline}).

    \begin{thm}\label{thm:Eaton}
        \textup{[\citet{Eaton1962-ai}]}. Let $\mathrm{RP1}$, $\mathrm{RT1}$, $\mathbf{p} \in \mathcal{S}^*$, $\mathbf{p} \neq \boldsymbol{0}$,
        \begin{equation*}
            \boldsymbol{q}(t, \mathbf{p}) \overset{\mathrm{def}}{=} \left(\frac{\boldsymbol{s}_{\mathcal{G}(t)}(-\mathbf{p}) - \boldsymbol{s}_{\mathcal{R}(t)}(\mathbf{p})}{\lVert\boldsymbol{s}_{\mathcal{G}(t)}(-\mathbf{p}) - \boldsymbol{s}_{\mathcal{R}(t)}(\mathbf{p})\rVert}\right)^\top - \mathbf{p},
        \end{equation*}
        $\rho_{\mathrm{lower}}(t, \mathbf{p}) \geq 0$, $\delta(t, \mathbf{p}) > 0$. Then there exists $\tilde{\gamma} \in \mathbb{R}^+$ such that for all $\gamma \in (0, \tilde{\gamma})$:
        \begin{equation}\label{eq:eaton_tmp}
            \rho_{\mathrm{lower}}(t, \mathbf{p} + \gamma\boldsymbol{q}(t, \mathbf{p})) > \rho_{\mathrm{lower}}(t, \mathbf{p}).
        \end{equation}
    \end{thm}
    \begin{pf}
        Lemma~\ref{lem:lower_diff_p} gives
        \begin{multline*}
            \nabla_{\mathbf{p}, \boldsymbol{q}(t, \mathbf{p})}\rho_{\mathrm{lower}}(t, \mathbf{p}) = \boldsymbol{q}(t, \mathbf{p})\frac{\partial\rho_{\mathrm{lower}}}{\partial\mathbf{p}}(t, \mathbf{p})\\
            = \boldsymbol{q}(t, \mathbf{p})\left(\frac{\boldsymbol{s}_{\mathcal{G}(t)}(-\mathbf{p}) - \boldsymbol{s}_{\mathcal{R}(t)}(\mathbf{p})}{\lVert\mathbf{p}\rVert} - \frac{\mathbf{p}^\top}{\lVert\mathbf{p}\rVert^2}\rho_{\mathrm{lower}}(t, \mathbf{p})\right)\\
            = \frac{\rho_{\mathrm{upper}}(t, \mathbf{p})}{\lVert\mathbf{p}\rVert} - \frac{\rho_{\mathrm{lower}}^2(t, \mathbf{p})}{\lVert\mathbf{p}\rVert\rho_{\mathrm{upper}}(t, \mathbf{p})}\\
            = \frac{\delta(t, \mathbf{p})(\rho_{\mathrm{upper}}(t, \mathbf{p}) + \rho_{\mathrm{lower}}(t, \mathbf{p}))}{\lVert\mathbf{p}\rVert\rho_{\mathrm{upper}}(t, \mathbf{p})} > 0.
        \end{multline*}
        The last inequality is true because $\delta(t, \mathbf{p}) > 0$ and $\rho_{\mathrm{upper}}(t, \mathbf{p}) > \rho_{\mathrm{lower}}(t, \mathbf{p}) \geq 0$. The directional derivative is continuous and positive. Thus, there exists $\tilde{\gamma} \in \mathbb{R}^+$ such that for all $\gamma \in [0, \tilde{\gamma})$ the directional derivative is positive. Hence,~\eqref{eq:eaton_tmp} holds for all $\gamma \in [0, \tilde{\gamma})$.
    \end{pf}
    
    The following property of $\boldsymbol{q}(t, \mathbf{p})$ will be useful to us. Let $\gamma \in \mathbb{R}$, $\lVert\mathbf{p}\rVert = 1$. Then
    \begin{multline}\label{eq:til_p_norm}
         \lVert\mathbf{p} + \gamma\boldsymbol{q}(t, \mathbf{p})\rVert^2 = \lVert\mathbf{p}\rVert^2 + 2\gamma\boldsymbol{q}(t, \mathbf{p})\mathbf{p}^\top + \gamma^2\lVert\boldsymbol{q}(t, \mathbf{p})\rVert^2\\
            = 1 + 2\gamma\left(\frac{\rho_{\mathrm{lower}}(t, \mathbf{p})}{\rho_{\mathrm{upper}}(t, \mathbf{p})} - 1\right) + 2\gamma^2\left(1 - \frac{\rho_{\mathrm{lower}}(t, \mathbf{p})}{\rho_{\mathrm{upper}}(t, \mathbf{p})}\right)\\
            = 1 - 2\gamma(1 - \gamma)\left(1 - \frac{\rho_{\mathrm{lower}}(t, \mathbf{p})}{\rho_{\mathrm{upper}}(t, \mathbf{p})}\right)\\
            = 1 - 2\gamma(1 - \gamma)\frac{\delta(t, \mathbf{p})}{\rho_{\mathrm{upper}}(t, \mathbf{p})}.
     \end{multline}

    \subsection{Time boosting}

\begin{figure}
        \begin{center}
            \includegraphics[width=0.4\textwidth]{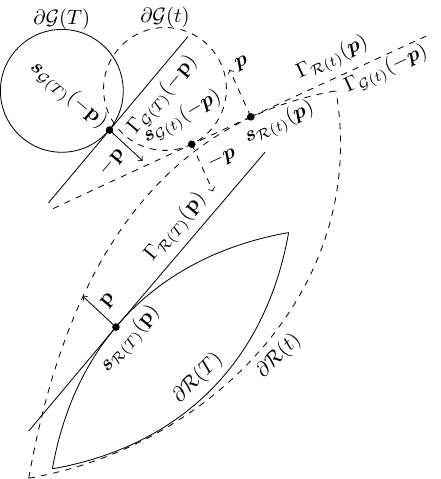}
            \caption{Evaluating of the boosting-time function $t = F(T, \mathbf{p})$. Here, $\boldsymbol{p} = \boldsymbol{p}(t; T, \mathbf{p})$}
            \label{fig:F}
        \end{center}
    \end{figure}

    The following function plays an important role for the algorithms:
    \begin{equation*}
        F(T, \mathbf{p}) \overset{\mathrm{def}}{=} \min\{t \in [T, +\infty):\: \rho_{\mathrm{lower}}(t, \boldsymbol{p}(t; T, \mathbf{p})) \leq 0\}.
    \end{equation*}
    We will call it {\it the boosting-time function}. If $\rho_{\mathrm{lower}}(T, \mathbf{p}) \geq 0$, i.e. $\mathcal{R}(t)$ and $\mathcal{G}(t)$ are separated, then $t = F(T, \mathbf{p}) \geq T$ is a minimal time moment such that
    \begin{equation*}
        \Gamma_{\mathcal{R}(t)}(\boldsymbol{p}(t; T, \mathbf{p})) = \Gamma_{\mathcal{G}(t)}(-\boldsymbol{p}(t; T, \mathbf{p})).   
    \end{equation*}
    Thus, $\mathcal{R}(t)$ and $\mathcal{G}(t)$ share the tangent hyperplane $\Gamma_{\mathcal{R}(t)}(\boldsymbol{p}(t; T, \mathbf{p}))$ (see~Fig.~\ref{fig:F}) and
    \begin{equation}\label{eq:F_main_prop}
        \rho_{\mathrm{lower}}(F(T, \mathbf{p}), \boldsymbol{p}(F(T, \mathbf{p}); T, \mathbf{p})) = 0.
    \end{equation}
    We conclude, that $\mathcal{R}(T)$ and $\mathcal{G}(T)$ are separated at $T$ and this property persists later at $t = F(T, \mathbf{p})$. This property of the function was first discovered and described by \citet{Neustadt1960-ki}.
        
    The ways of calculating $F(T, \mathbf{p})$ and the encountered difficulties will be described later. Now we note that the direct calculation of $F(T, \mathbf{p})$ is difficult without numerical integration, even if the boundary of the reachable set is described analytically. Further we describe a simpler method for time boosting that persists the separability of $\mathcal{R}(t)$ and $\mathcal{G}(t)$.

    The function given by
    \begin{equation*}
        f(T, \mathbf{p}) = 
        \begin{cases}
            T + \frac{\rho_{\mathrm{lower}}(T, \mathbf{p})}{v_{\mathcal{G}} + v_{\mathcal{R}}}, \quad &\rho_{\mathrm{lower}}(T, \mathbf{p}) \geq 0;\\
            T, \quad &\rho_{\mathrm{lower}}(T, \mathbf{p}) < 0
        \end{cases}
    \end{equation*}
    is {\it the simple boosting-time function}.
    
    \begin{lem}\label{lem:f_separ}
        Let $\mathrm{RP1}$, $\mathrm{RP2}$, $\mathrm{RT4}$, $\mathbf{p} \in \mathcal{S}^*$, $\mathbf{p} \neq \boldsymbol{0}$, $\rho_{\mathrm{lower}}(T, \mathbf{p}) > 0$. Then
        \begin{equation*}
            \rho_{\mathrm{lower}}(t, \mathbf{p}) > 0, \quad t \in [T, f(T, \mathbf{p})),
        \end{equation*}
        and $\rho_{\mathrm{lower}}(f(T, \mathbf{p}), \mathbf{p}) \geq 0$.
    \end{lem}
    \begin{pf}
        Using RP2, RT4 and Lemma~\ref{lem:lower_diff_t}, we deduce
        \begin{multline*}
            \left|\frac{\partial\rho_{\mathrm{lower}}}{\partial t}(\tau, \mathbf{p})\right|\\
            = \left|\frac{\mathbf{p}}{\lVert\mathbf{p}\rVert}(\boldsymbol{v}_{\mathcal{G}(\tau)}(-\mathbf{p}) - \boldsymbol{A}(\tau)\boldsymbol{s}_{\mathcal{R}(\tau)}(\mathbf{p}) - \boldsymbol{u}_E(\mathbf{p}))\right|\\
            \leq \lVert\boldsymbol{v}_{\mathcal{G}(\tau)}(-\mathbf{p}) - \boldsymbol{A}(\tau)\boldsymbol{s}_{\mathcal{R}(\tau)}(\mathbf{p}) - \boldsymbol{u}_E(\mathbf{p})\rVert\\
            \leq \lVert\boldsymbol{v}_{\mathcal{G}(\tau)}(-\mathbf{p})\rVert + \lVert\boldsymbol{A}(\tau)\boldsymbol{s}_{\mathcal{R}(\tau)}(\mathbf{p}) + \boldsymbol{u}_E(\mathbf{p})\rVert\\
            \leq v_{\mathcal{G}} + v_{\mathcal{R}}.
        \end{multline*}
        Thus,
        \begin{multline*}
            \rho_{\mathrm{lower}}(t, \mathbf{p}) = \rho_{\mathrm{lower}}(T, \mathbf{p}) + \int\limits_{T}^{t}\frac{\partial\rho_{\mathrm{lower}}}{\partial t}(\tau, \mathbf{p})\mathrm{d}\tau\\
            \geq \rho_{\mathrm{lower}}(T, \mathbf{p}) - \int\limits_{T}^{t}\left|\frac{\partial\rho_{\mathrm{lower}}}{\partial t}(\tau, \mathbf{p})\right|\mathrm{d}\tau\\
            \geq \rho_{\mathrm{lower}}(T, \mathbf{p}) - (t - T)(v_{\mathcal{G}} + v_{\mathcal{R}}).\\
        \end{multline*}
        It implies $\rho_{\mathrm{lower}}(t, \mathbf{p}) \geq 0$ for $t = f(T, \mathbf{p})$ and $\rho_{\mathrm{lower}}(t, \mathbf{p}) > 0$ for $t \in [T, f(T, \mathbf{p}))$.
    \end{pf}



    \section{Algorithms of computing the distance}\label{sec:MDP}

    Each algorithm for solving the general problem (MTPLS) uses the initial value of the support vector $\mathbf{p}_0 \in \mathcal{S}^*$ ($\mathbf{p}_0 \neq \boldsymbol{0}$) as input data. For some algorithms, this vector must be chosen so that the hyperplane $\Gamma_{\mathcal{G}(0)}(-\mathbf{p}_0)$ separates the initial state $\mathbf{s}_0$ from the target set $\mathcal{G}(0)$, i.e. $\rho_{\mathrm{lower}}(0, \mathbf{p}_0) \geq 0$. Note, that if $\delta(0, \mathbf{p}_0) = 0$, then $\rho_{\mathrm{lower}}(0, \mathbf{p}_0) = \rho(0) \geq 0$. Therefore, solving the problem of finding the closest point $\boldsymbol{s}_{\mathcal{G}(0)}(-\mathbf{p})$ on $\mathcal{G}(0)$ to $\mathbf{s}_0$ by tuning $\mathbf{p}$ becomes actual. A more general case of this problem is the MDP for $\mathcal{R}(t)$ and $\mathcal{G}(t)$.

    \subsection{Gilbert-Johnson-Keerthi distance algorithm}

    The most famous algorithm for finding the distance between two convex sets is GJK distance algorithm \citep{Gilbert1988-ug, Gilbert1990-wx}. This algorithm actively uses the contact functions\footnote{The original GJK algorithm also applies to non-strictly convex sets, i.e. $\boldsymbol{s}_{\mathcal{R}(t)}$, $\boldsymbol{s}_{\mathcal{G}(t)}$ are not single-valued functions and return support sets instead of support points. For convenience, we will assume that all sets are strictly convex.} $\boldsymbol{s}_{\mathcal{R}(t)}$, $\boldsymbol{s}_{\mathcal{G}(t)}$ and Johnson's algorithm. 

    Let $\mathcal{V}$ be a finite set of states and $\mathrm{card}(\mathcal{V}) \leq n + 1$. Johnson's algorithm computes the distance from the convex hull of $\mathcal{V}$ to the origin $\boldsymbol{0}$. The output of Johnson's algorithm $(\mathbf{s}, \tilde{\mathcal{V}}) = \mathrm{Johnson}(\mathcal{V})$ is a tuple of two elements. The first one is the nearest point $\mathbf{s}$ on $\mathrm{conv}(\mathcal{V})$ to $\boldsymbol{0}$. The second one $\tilde{\mathcal{V}}$ is a subset of $\mathcal{V}$ such that $\mathbf{s} \in \mathrm{conv}(\tilde{\mathcal{V}})$. If $\boldsymbol{0} \notin \mathrm{conv}(\mathcal{V})$, then $\mathrm{card}(\tilde{\mathcal{V}}) \leq n$.
    
    \begin{algorithm}
    	\caption{GJK distance algorithm: $\mathbf{GJK}_{\varepsilon}(t, \mathbf{p})$}
    	\label{alg:GJK}
    	\textbf{Input:} $\varepsilon \in \mathbb{R}^+$, $t \in \mathbb{R}^+_0$, $\mathbf{p} \in \mathcal{S}^*$\\
        \textbf{Require:} $\mathrm{RP1}$, $\mathrm{RT1}$, $\mathbf{p} \neq \boldsymbol{0}$
    	\begin{algorithmic}[1]
            \State $\mathbf{s} \gets \boldsymbol{s}_{\mathcal{R}(t)}(\mathbf{p}) - \boldsymbol{s}_{\mathcal{G}(t)}(-\mathbf{p})$
            \State $\mathcal{V} \gets \{\mathbf{s}\}$
            \While{$\lVert\mathbf{s}\rVert > \varepsilon$ \textbf{and} $\lVert\mathbf{s}\rVert - \rho_{\mathrm{lower}}(t, -\mathbf{s}^\top) > \varepsilon$}
                \State $\mathcal{V} \gets \mathcal{V} \cup \{\boldsymbol{s}_{\mathcal{R}(t)}(-\mathbf{s}^\top) - \boldsymbol{s}_{\mathcal{G}(t)}(\mathbf{s}^\top)\}$
                \State $\mathbf{s}, \mathcal{V} \gets \mathrm{Johnson}(\mathcal{V})$
            \EndWhile
    	\end{algorithmic}	
    	\textbf{Output:} $\mathbf{s}$
    \end{algorithm}
    
    The GJK algorithm works as follows. The original problem of finding the distance between $\mathcal{G}(t)$ and $\mathcal{R}(t)$ is replaced by the problem of finding the distance between $\mathcal{P}(t) = \mathcal{R}(t) - \mathcal{G}(t)$ and $\boldsymbol{0}$. At each iteration, there is a set of vertices $\mathcal{V}$ of a convex polyhedron. Each vertex is a previously calculated point from the boundary of $\mathcal{P}(t)$. Using Johnson's algorithm, we compute the distance vector between the convex polyhedron and the origin. Johnson's algorithm also identifies the face or edge of the convex polyhedron on which the point closest to the origin lies. The new set of vertices $\mathcal{V}$ now consists of the vertices of this face/edge plus a new point from the boundary of $\mathcal{P}(t)$, for which the computed distance vector is the support vector.

    \begin{thm}\label{thm:gilbert-barr}
        \textup{[\citet{Gilbert1990-wx}]}. Let $\mathrm{RP1}$, $\mathrm{RT1}$, $\varepsilon \in \mathbb{R}^+$, $t \in \mathbb{R}^+_0$, $\mathbf{p} \in \mathcal{S}^*$ ($\mathbf{p} \neq \boldsymbol{0}$). Then, Algorithm~\ref{alg:GJK} terminates in a finite number of iterations. Let $\mathbf{GJK}_{\varepsilon}(t, \mathbf{p}) \in \mathcal{S}$ be an output of Algorithm~\ref{alg:GJK}. Then $\lVert\mathbf{GJK}_{\varepsilon}(t, \mathbf{p})\rVert - \rho(t) \leq \varepsilon$. Let
        \begin{equation*}
            \tilde{\mathbf{s}} = \lim_{\varepsilon \to +0} \mathbf{GJK}_{\varepsilon}(t, \mathbf{p}).
        \end{equation*}
        If $\tilde{\mathbf{s}} \neq \boldsymbol{0}$, then $\delta(t, -\tilde{\mathbf{s}}^\top) = 0$. If $\tilde{\mathbf{s}} = \boldsymbol{0}$, then $\mathcal{R}(t) \cap \mathcal{G}(t) \neq \varnothing$.
    \end{thm}
    \begin{pf}
        It follows directly from Theorem 1 of~\citet{Barr1969-ca} and \citet{Gilbert1990-wx}.
    \end{pf}
    \begin{cor}\label{cor:gjk_sep}
        Let $\mathbf{s} = \mathbf{GJK}_{\varepsilon}(t, \mathbf{p})$. If $\lVert\mathbf{s}\rVert > \varepsilon$, then $\rho_{\mathrm{lower}}(t, -\mathbf{s}^\top) > 0$.
    \end{cor}
    \begin{pf}
        If $\lVert\mathbf{s}\rVert > \varepsilon$, then the while-loop of Algorithm~\ref{alg:GJK} terminates due to
        \begin{equation*}
           \lVert\mathbf{s}\rVert - \rho_{\mathrm{lower}}(t, -\mathbf{s}^\top) \leq \varepsilon.
        \end{equation*}
        Using $\lVert\mathbf{s}\rVert > \varepsilon$, we obtain $\rho_{\mathrm{lower}}(t, -\mathbf{s}^\top) > 0$.
    \end{pf}

    If $\mathcal{R}(t) \cap \mathcal{G}(t) = \varnothing$, then the variable $\mathbf{s}$ of Algorithm~\ref{alg:GJK} is non-zero at any iteration, because $\mathbf{s} \in \mathcal{R}(t) - \mathcal{G}(t) \not\owns \boldsymbol{0}$. Algorithm~\ref{alg:GJKs} is a modified GJK algorithm, that uses $\mathbf{s} \neq \boldsymbol{0}$ and $\mathcal{R}(t) \cap \mathcal{G}(t) = \varnothing$. In contrast, it returns the element of conjugate space $\mathcal{S}^*$ instead of $\mathcal{S}$. We denote the output of Algorithm~\ref{alg:GJKs} by $\mathbf{GJK}^*_{\alpha}(t, \mathbf{p}) \in \mathcal{S}^*$. It will be demonstrated that the interface of Algorithm~\ref{alg:GJKs} is more convenient for the MTPLS.

    \begin{algorithm}
    	\caption{modified GJK: $\mathbf{GJK}^*_{\alpha}(t, \mathbf{p})$}
    	\label{alg:GJKs}
    	\textbf{Input:} $\alpha \in \mathbb{R}^+$, $t \in \mathbb{R}^+_0$, $\mathbf{p} \in \mathcal{S}^*$\\
        \textbf{Require:} $\mathrm{RP1}$, $\mathrm{RT1}$, $\mathbf{p} \neq \boldsymbol{0}$, $\mathcal{R}(t) \cap \mathcal{G}(t) = \varnothing$
    	\begin{algorithmic}[1]
            \State $\mathcal{V} \gets \{\boldsymbol{s}_{\mathcal{R}(t)}(\mathbf{p}) - \boldsymbol{s}_{\mathcal{G}(t)}(-\mathbf{p})\}$
            \While{$\delta(t, \mathbf{p}) > \alpha\rho_{\mathrm{lower}}(t, \mathbf{p})$}
                \State $\mathbf{s}, \mathcal{V} \gets \mathrm{Johnson}(\mathcal{V})$
                \State $\mathbf{p} \gets -\mathbf{s}^\top/\lVert\mathbf{s}\rVert$
                \State $\mathcal{V} \gets \mathcal{V} \cup \{\boldsymbol{s}_{\mathcal{R}(t)}(\mathbf{p}) - \boldsymbol{s}_{\mathcal{G}(t)}(-\mathbf{p})\}$
            \EndWhile
    	\end{algorithmic}	
    	\textbf{Output:} $\mathbf{p}$
    \end{algorithm}

    \begin{thm}
        Let $\mathrm{RP1}$, $\mathrm{RT1}$, $\alpha \in \mathbb{R}^+$, $t \in \mathbb{R}^+_0$, $\mathbf{p} \in \mathcal{S}^*$ ($\mathbf{p} \neq \boldsymbol{0}$), and $\mathcal{R}(t) \cap \mathcal{G}(t) = \varnothing$. Then, Algorithm~\ref{alg:GJKs} terminates in a finite number of iterations and $\delta(t, \hat{\mathbf{p}}) \leq \alpha\rho_{\mathrm{lower}}(t, \hat{\mathbf{p}})$, where $\hat{\mathbf{p}} = \mathbf{GJK}^*_{\alpha}(t, \mathbf{p})$. Let
        \begin{equation*}
            \tilde{\mathbf{p}} = \lim_{\alpha \to +0} \mathbf{GJK}^*_{\alpha}(t, \mathbf{p}).
        \end{equation*}
        Then, $\delta(t, \tilde{\mathbf{p}}) = 0$.
    \end{thm}
    \begin{pf}
        It follows directly from Theorem~\ref{thm:gilbert-barr} and $\lVert\mathbf{s}\rVert > 0$.
    \end{pf}

    \begin{cor}\label{cor:gjks_sep}
        $\rho_{\mathrm{lower}}(t, \tilde{\mathbf{p}}) > 0$ for $\tilde{\mathbf{p}} = \mathbf{GJK}^*_{\alpha}(t, \mathbf{p})$.
    \end{cor}
    \begin{pf}
        Algorithm~\ref{alg:GJKs} terminates with
        \begin{equation*}
            \delta(t, \tilde{\mathbf{p}}) \leq \alpha\rho_{\mathrm{lower}}(t, \tilde{\mathbf{p}}).
        \end{equation*}
        Using~$\mathcal{R}(t) \cap \mathcal{G}(t) = \varnothing$, we obtain $\rho_{\mathrm{upper}}(t, \tilde{\mathbf{p}}) \geq \rho(t) > 0$. Thus,
        \begin{equation*}
            (1 + \alpha)\rho_{\mathrm{lower}}(t, \tilde{\mathbf{p}}) \geq \rho_{\mathrm{upper}}(t, \tilde{\mathbf{p}}) > 0.
        \end{equation*}
    \end{pf}

    \subsection{Gilbert distance algorithm}

    Gilbert's distance algorithm is the predecessor of the GJK distance algorithm. Instead of a convex polyhedron with vertices $\mathcal{V}$ it uses only a line segment or a single point. We denote the output of Algorithm~\ref{alg:G} by $\mathbf{G}_{\alpha}(t, \mathbf{p}) \in \mathcal{S}^*$
    
    \begin{algorithm}
    	\caption{Gilbert distance algorithm: $\mathbf{G}_{\alpha}(t, \mathbf{p})$}
    	\label{alg:G}
    	\textbf{Input:} $\alpha \in \mathbb{R}^+$, $t \in \mathbb{R}^+_0$, $\mathbf{p} \in \mathcal{S}^*$\\
        \textbf{Require:} $\mathrm{RP1}$, $\mathrm{RT1}$, $\mathbf{p} \neq \boldsymbol{0}$, $\mathcal{R}(t) \cap \mathcal{G}(t) = \varnothing$
    	\begin{algorithmic}[1]
                \State $\mathbf{s} \gets \boldsymbol{s}_{\mathcal{R}(t)}(\mathbf{p}) - \boldsymbol{s}_{\mathcal{G}(t)}(-\mathbf{p})$
                \While{$\delta(t, \mathbf{p}) > \alpha\rho_{\mathrm{lower}}(t, \mathbf{p})$}
                    \State $\mathbf{p} \gets -\mathbf{s}^\top/\lVert\mathbf{s}\rVert$
                    \State $\tilde{\mathbf{s}} \gets \mathbf{s} - (\boldsymbol{s}_{\mathcal{R}(t)}(\mathbf{p}) - \boldsymbol{s}_{\mathcal{G}(t)}(-\mathbf{p}))$
                    \State $\mathbf{s} \gets \mathbf{s} + \tilde{\mathbf{s}}\min\left(1, \frac{\mathbf{s}^\top\tilde{\mathbf{s}}}{\tilde{\mathbf{s}}^2}\right)$
                \EndWhile
    	\end{algorithmic}	
    	\textbf{Output:} $\mathbf{p}$
    \end{algorithm}
    
    \begin{thm}
        \textup{[\citet{Gilbert1966-sc}]} Let $\mathrm{RP1}$, $\mathrm{RT1}$, $\alpha \in \mathbb{R}^+$, $t \in \mathbb{R}^+_0$, $\mathbf{p} \in \mathcal{S}^*$ ($\mathbf{p} \neq \boldsymbol{0}$), and $\mathcal{R}(t) \cap \mathcal{G}(t) = \varnothing$. Then, Algorithm~\ref{alg:G} terminates in a finite number of iterations and $\delta(t, \hat{\mathbf{p}}) \leq \alpha\rho_{\mathrm{lower}}(t, \hat{\mathbf{p}})$, where $\hat{\mathbf{p}} = \mathbf{G}_{\alpha}(t, \mathbf{p})$. Let
        \begin{equation*}
            \tilde{\mathbf{p}} = \lim_{\alpha \to +0} \mathbf{G}_{\alpha}(t, \mathbf{p}).
        \end{equation*}
        Then, $\delta(t, \tilde{\mathbf{p}}) = 0$.
    \end{thm}
    \begin{pf}
        It follows directly from Theorem~2 of \cite{Gilbert1966-sc}.
    \end{pf}

    \begin{cor}\label{cor:g_sep}
        $\rho_{\mathrm{lower}}(t, \tilde{\mathbf{p}}) > 0$ for $\tilde{\mathbf{p}} = \mathbf{G}_{\alpha}(t, \mathbf{p})$.
    \end{cor}
    \begin{pf}
        Algorithm~\ref{alg:G} terminates with
        \begin{equation*}
            \delta(t, \tilde{\mathbf{p}}) \leq \alpha\rho_{\mathrm{lower}}(t, \tilde{\mathbf{p}}).
        \end{equation*}
        Using~$\mathcal{R}(t) \cap \mathcal{G}(t) = \varnothing$, we obtain $\rho_{\mathrm{upper}}(t, \tilde{\mathbf{p}}) \geq \rho(t) > 0$. Thus,
        \begin{equation*}
            (1 + \alpha)\rho_{\mathrm{lower}}(t, \tilde{\mathbf{p}}) \geq \rho_{\mathrm{upper}}(t, \tilde{\mathbf{p}}) > 0.
        \end{equation*}
    \end{pf}
    
    \subsection{Steepest ascent distance algorithm}

    We now describe an original algorithm for computing the distance that is simpler than Algorithm~\ref{alg:GJKs}, since it doesn't use Johnson's algorithm. This algorithm is based on maximizing $\rho_{\mathrm{lower}}(t, \cdot)$ by steepest ascent method. Theorem~\ref{thm:Eaton} states that $\boldsymbol{q}(t, \mathbf{p})$ is the direction of ascending for $\rho_{\mathrm{lower}}(t, \cdot)$ at $\mathbf{p}$. Next, we find out the conditions of convergence of the steepest ascent algorithm and set the step-size guaranteeing convergence. The last one is a novel result of the paper.
    
    First we obtain several auxiliary results. Let
    \begin{multline*}
        \xi(\gamma; t, \mathbf{p}) \overset{\mathrm{def}}{=} \lVert\mathbf{p} + \gamma\boldsymbol{q}(t, \mathbf{p})\rVert\rho_{\mathrm{lower}}(t, \mathbf{p} + \gamma\boldsymbol{q}(t, \mathbf{p})) \\
        - \rho_{\mathrm{lower}}(t, \mathbf{p})
    \end{multline*}
     
    \begin{lem}\label{lem:xi_prop}
       Let $\mathrm{RP1}$, $\mathrm{RT1}$, $t \in \mathbb{R}^+_0$, $\mathbf{p} \in \mathcal{S}^*$, $\lVert\mathbf{p}\rVert = 1$, $\delta(t, \mathbf{p}) > 0$, and $\rho_{\mathrm{lower}}(t, \mathbf{p}) \geq 0$. Then
       \begin{multline*}
           \xi'(\gamma; t, \mathbf{p}) \overset{\mathrm{def}}{=} \frac{\partial\xi}{\partial\gamma}(\gamma; t, \mathbf{p})\\
           = \boldsymbol{q}(t, \mathbf{p})(\boldsymbol{s}_{\mathcal{G}(t)}(-(\mathbf{p} + \gamma\boldsymbol{q}(t, \mathbf{p}))) - \boldsymbol{s}_{\mathcal{R}(t)}(\mathbf{p} + \gamma\boldsymbol{q}(t, \mathbf{p}))),
       \end{multline*}
       $\xi'(0; t, \mathbf{p}) > 0$, and $\xi'(\cdot; t, \mathbf{p})$ is a non-increasing function.
    \end{lem}
    \begin{pf}
       Let denote $\tilde{\mathbf{p}} = \mathbf{p} + \gamma\boldsymbol{q}(t, \mathbf{p})$. Using~\eqref{eq:til_p_norm}, we obtain
       \begin{multline*}
           \lVert\tilde{\mathbf{p}}\rVert^2 = 1 - 2\gamma(1 - \gamma)\frac{\delta(t, \mathbf{p})}{\rho_{\mathrm{upper}}(t, \mathbf{p})}\\
           = 1 - (1 - \gamma^2 - (1 - \gamma)^2)\frac{\delta(t, \mathbf{p})}{\rho_{\mathrm{upper}}(t, \mathbf{p})}\\
           > 1 - \frac{\delta(t, \mathbf{p})}{\rho_{\mathrm{upper}}(t, \mathbf{p})} = \frac{\rho_{\mathrm{lower}}(t, \mathbf{p})}{\rho_{\mathrm{upper}}(t, \mathbf{p})} \geq 0,
       \end{multline*}
       since $\rho_{\mathrm{lower}}(t, \mathbf{p}) \geq 0$ and $\delta(t, \mathbf{p}) > 0$. Thus, $\tilde{\mathbf{p}} \neq \boldsymbol{0}$. Using~\eqref{eq:xi_will}, we obtain
       \begin{multline*}
            \xi'(\gamma; t, \mathbf{p}) = \nabla_{\tilde{\mathbf{p}}, \boldsymbol{q}(t, \mathbf{p})}(\lVert\tilde{\mathbf{p}}\rVert\rho_{\mathrm{lower}}(t, \tilde{\mathbf{p}}))\\
            = \boldsymbol{q}(t, \mathbf{p})\frac{\partial}{\partial\tilde{\mathbf{p}}}(\lVert\tilde{\mathbf{p}}\rVert\rho_{\mathrm{lower}}(t, \tilde{\mathbf{p}}))\\
            = \boldsymbol{q}(t, \mathbf{p})(\boldsymbol{s}_{\mathcal{G}(t)}(-(\mathbf{p} + \gamma\boldsymbol{q}(t, \mathbf{p}))) - \boldsymbol{s}_{\mathcal{R}(t)}(\mathbf{p} + \gamma\boldsymbol{q}(t, \mathbf{p}))).
       \end{multline*}
       Now we check that $\xi'(0; t, \mathbf{p}) > 0$:
       \begin{multline*}
           \xi'(0; t, \mathbf{p}) = \boldsymbol{q}(t, \mathbf{p})(\boldsymbol{s}_{\mathcal{G}(t)}(-\mathbf{p}) - \boldsymbol{s}_{\mathcal{R}(t)}(\mathbf{p}))\\
           = \rho_{\mathrm{upper}}(t, \mathbf{p}) - \rho_{\mathrm{lower}}(t, \mathbf{p}) = \delta(t, \mathbf{p}) > 0.
       \end{multline*}
       To prove non-increasing of $\xi'(\cdot; t, \mathbf{p})$ we consider the following difference:
       \begin{multline*}
           \xi'(\gamma + b; t, \mathbf{p}) - \xi'(\gamma; t, \mathbf{p}) \\
           = \boldsymbol{q}(t, \mathbf{p})(\boldsymbol{s}_{\mathcal{G}(t)}(-(\tilde{\mathbf{p}} + b\boldsymbol{q}(t, \mathbf{p}))) - \boldsymbol{s}_{\mathcal{G}(t)}(-\tilde{\mathbf{p}}))\\
           + \boldsymbol{q}(t, \mathbf{p})(\boldsymbol{s}_{\mathcal{R}(t)}(\tilde{\mathbf{p}}) - \boldsymbol{s}_{\mathcal{R}(t)}(\tilde{\mathbf{p}} + b\boldsymbol{q}(t, \mathbf{p})))\\
           = \frac{1}{b}(-(\tilde{\mathbf{p}} + b\boldsymbol{q}(t, \mathbf{p})))(\boldsymbol{s}_{\mathcal{G}(t)}(-\tilde{\mathbf{p}}) - \boldsymbol{s}_{\mathcal{G}(t)}(-(\tilde{\mathbf{p}} + b\boldsymbol{q}(t, \mathbf{p}))))\\
           + \frac{1}{b}(-\tilde{\mathbf{p}})(\boldsymbol{s}_{\mathcal{G}(t)}(-(\tilde{\mathbf{p}} + b\boldsymbol{q}(t, \mathbf{p}))) - \boldsymbol{s}_{\mathcal{G}(t)}(-\tilde{\mathbf{p}}))\\
           + \frac{1}{b}(\tilde{\mathbf{p}} + b\boldsymbol{q}(t, \mathbf{p}))(\boldsymbol{s}_{\mathcal{R}(t)}(\tilde{\mathbf{p}}) - \boldsymbol{s}_{\mathcal{R}(t)}(\tilde{\mathbf{p}} + b\boldsymbol{q}(t, \mathbf{p})))\\
           + \frac{1}{b}\tilde{\mathbf{p}}(\boldsymbol{s}_{\mathcal{R}(t)}(\tilde{\mathbf{p}} + b\boldsymbol{q}(t, \mathbf{p})) - \boldsymbol{s}_{\mathcal{R}(t)}(\tilde{\mathbf{p}}))
       \end{multline*}
       Using~\eqref{eq:conv}, we conclude that if $b$ is positive, then the difference is non-positive.
    \end{pf}

    \begin{lem}\label{lem:supp_vect_alig_finit}
        Let $\mathrm{RP1}$, $\mathrm{RT1}$, $\gamma \in \mathbb{R}^+$, $t \in \mathbb{R}^+_0$, $\mathbf{p} \in \mathcal{S}^*$, $\lVert\mathbf{p}\rVert = 1$, $\delta(t, \mathbf{p}) > 0$, and $\rho_{\mathrm{lower}}(t, \mathbf{p}) \geq 0$. There exists minimal $N \in \mathbb{N}$ such that for any $n \in \mathbb{N}$ if $n \geq N$, then $\xi'(\gamma/2^{n - 1}; t, \mathbf{p}) \geq 0$, but if $n < N$, then $\xi'(\gamma/2^{n - 1}; t, \mathbf{p}) < 0$.
    \end{lem}
    \begin{pf}
        Lemmas~\ref{lem:cont_cont},~\ref{lem:xi_prop} confirm that $\xi'(\cdot; t, \mathbf{p})$ is a continuous function. Lemma~\ref{lem:xi_prop} also states that $\xi'(\cdot; t, \mathbf{p})$ is a non-increasing function and $\xi'(0; t, \mathbf{p}) > 0$. Thus, either $\xi'(\tilde{\gamma}; t, \mathbf{p})$ is positive for any $\tilde{\gamma} \in \mathbb{R}^+$ or there exists $\gamma_{\min} \in \mathbb{R}^+$ such that $\xi'(\tilde{\gamma}; t, \mathbf{p}) \geq 0$ if $\tilde{\gamma} \in [0, \gamma_{\min}]$, and $\xi'(\tilde{\gamma}; t, \mathbf{p}) < 0$ if $\tilde{\gamma} > \gamma_{\min}$. In the first case, we set $N = 1$, because $\xi'(\gamma/2^{n - 1}; t, \mathbf{p}) > 0$ for all $n \geq N$. In the second case, we use the minimal $N \in \mathbb{N}$ such that $\gamma/2^{N - 1} \leq \gamma_{\min}$.
    \end{pf}
    
    Lemma~\ref{lem:supp_vect_alig_finit} guarantees that
    \begin{equation}\label{eq:N_def}
        N(t, \mathbf{p}) \overset{\mathrm{def}}{=} \min\{n \in \mathbb{N}:\: \xi'(\gamma/2^{n - 1}; t, \mathbf{p}) \geq 0\}
    \end{equation}
    can be computed iteratively with a finite number of iterations.

    \begin{lem}\label{lem:rho_lower_incr}
       Let $\mathrm{RP1}$, $\mathrm{RT1}$, $\gamma \in (0, 1]$, $t \in \mathbb{R}^+_0$, $\mathbf{p} \in \mathcal{S}^*$, $\lVert\mathbf{p}\rVert = 1$, $\delta(t, \mathbf{p}) > 0$, and $\rho_{\mathrm{lower}}(t, \mathbf{p}) \geq 0$. Then
       \begin{equation*}
           \rho_{\mathrm{lower}}(t, \mathbf{p} + \gamma/2^{N(t, \mathbf{p}) - 1}\boldsymbol{q}(t, \mathbf{p})) > \rho_{\mathrm{lower}}(t, \mathbf{p}).
       \end{equation*}
    \end{lem}
    \begin{pf}
        By the definition $\xi'(\gamma/2^{N(t, \mathbf{p}) - 1}; t, \mathbf{p}) \geq 0$ and $\xi'(\tilde{\gamma}; t, \mathbf{p}) \geq 0$ for all $\tilde{\gamma} \in [0, \gamma/2^{N(t, \mathbf{p}) - 1}]$. Using $\xi(0; t, \mathbf{p}) = 0$ and $\xi'(0; t, \mathbf{p}) > 0$, we conclude that $\xi(\gamma/2^{N(t, \mathbf{p}) - 1}; t, \mathbf{p}) > 0$. Hence
        \begin{multline*}
            \lVert\mathbf{p} + \gamma/2^{N(t, \mathbf{p}) - 1}\boldsymbol{q}(t, \mathbf{p})\rVert\rho_{\mathrm{lower}}(t, \mathbf{p} + \gamma/2^{N(t, \mathbf{p}) - 1}\boldsymbol{q}(t, \mathbf{p})) \\
            > \rho_{\mathrm{lower}}(t, \mathbf{p}).
        \end{multline*}
        It remains to prove that $\lVert\mathbf{p} + \gamma/2^{N(t, \mathbf{p}) - 1}\boldsymbol{q}(t, \mathbf{p})\rVert \leq 1$. It follows from~\eqref{eq:til_p_norm}, since
        \begin{equation*}
            \frac{\delta(t, \mathbf{p})}{\rho_{\mathrm{upper}}(t, \mathbf{p})} > 0, \quad \gamma/2^{N(t, \mathbf{p}) - 1}(1 - \gamma/2^{N(t, \mathbf{p}) - 1}) \geq 0.
         \end{equation*}
    \end{pf}

    Lemma~\ref{lem:rho_lower_incr} is a constructive refinement of Theorem~\ref{thm:Eaton}. It declares an unambiguous way to compute the step-size that guarantees the increase of $\rho_{\mathrm{lower}}(t, \cdot)$. Condition~\eqref{eq:eaton_tmp} of Theorem~\ref{thm:Eaton} does not constructively define $\tilde{\gamma}$, since it is possible that~\eqref{eq:eaton_tmp} is satisfied for some $\hat{\gamma} = \gamma/2^{n - 1}$, but it does not follow that it is satisfied for all $\gamma \in (0, \hat{\gamma})$. This fact is the main reason why the Neustadt-Eaton algorithm (described in the next section) does not have a constructive convergence proof.

    Now we are ready to describe a novel algorithm for computing the distance function. Algorithm~\ref{alg:SA} is a steepest ascent algorithm that uses specific rules of tuning the step-size.
    \begin{algorithm}
    	\caption{Steepest ascent algorithm: $\mathbf{SA}_{\alpha}(t, \mathbf{p}; \gamma)$}
    	\label{alg:SA}
    	\textbf{Input:} $\alpha \in \mathbb{R}^+$, $t \in \mathbb{R}^+_0$, $\mathbf{p} \in \mathcal{S}^*$, $\gamma \in (0, 1)$\\
        \textbf{Require:} $\mathrm{RP1}$, $\mathrm{RT1}$, $\lVert\mathbf{p}\rVert = 1$, $\rho_{\mathrm{lower}}(t, \mathbf{p}) > 0$
    	\begin{algorithmic}[1]
            \While{$\delta(t, \mathbf{p}) > \alpha\rho_{\mathrm{lower}}(t, \mathbf{p})$}
                \While{$\xi'(\gamma; t, \mathbf{p}) < 0$}
                    \State $\gamma \gets \gamma / 2$
                \EndWhile
                \State $\mathbf{p} \gets \frac{\mathbf{p} + \gamma\boldsymbol{q}(t, \mathbf{p})}{\lVert\mathbf{p} + \gamma\boldsymbol{q}(t, \mathbf{p})\rVert}$
            \EndWhile
    	\end{algorithmic}	
    	\textbf{Output:} $\mathbf{p}$
    \end{algorithm}
    Lemma~\ref{lem:rho_lower_incr} states that the lower estimation ascents and persists positive at any stage of Algorithm~\ref{alg:SA}. Lemma~\ref{lem:supp_vect_alig_finit} guarantees the termination of the nested while-loop of Algorithm~\ref{alg:SA}.
    
    Now we prove that Algorithm~\ref{alg:SA} terminates in a finite number of iterations. We suppose that Algorithm~\ref{alg:SA} generate the following sequence of support vectors:
    \begin{equation*}
        \mathbf{p}_0 \overset{\mathrm{def}}{=} \mathbf{p},\quad \mathbf{p}_{m + 1} \overset{\mathrm{def}}{=} \frac{\mathbf{p}_m + \gamma/2^{N(t, \mathbf{p}_m) - 1}\boldsymbol{q}(t, \mathbf{p}_m)}{\lVert\mathbf{p}_m + \gamma/2^{N(t, \mathbf{p}_m) - 1}\boldsymbol{q}(t, \mathbf{p}_m)\rVert}.
    \end{equation*}
    We denote the output of Algorithm~\ref{alg:SA} by $\mathbf{SA}_{\alpha}(t, \mathbf{p}; \gamma)$.
   
    \begin{thm}\label{thm:dist_est_finit}
        Let $\mathrm{RP1}$, $\mathrm{RT1}$, $\alpha \in \mathbb{R}^+$, $t \in \mathbb{R}^+_0$, $\mathbf{p} \in \mathcal{S}^*$, $\gamma \in (0, 1)$, $\lVert\mathbf{p}\rVert = 1$, $\rho_{\mathrm{lower}}(t, \mathbf{p}) \geq 0$, and $\mathcal{R}(t) \cap \mathcal{G}(t) = \varnothing$. Then, Algorithm~\ref{alg:SA} terminates in a finite number of iterations and $\delta(t, \hat{\mathbf{p}}) \leq \alpha\rho_{\mathrm{lower}}(t, \hat{\mathbf{p}})$, where $\hat{\mathbf{p}} = \mathbf{SA}_\alpha(t, \mathbf{p}; \gamma)$. Let
        \begin{equation*}
            \tilde{\mathbf{p}} = \lim_{\alpha \to +0} \mathbf{SA}_\alpha(t, \mathbf{p}; \gamma).
        \end{equation*}
        Then, $\delta(t, \tilde{\mathbf{p}}) = 0$.
    \end{thm}
   \begin{pf}
       If $\delta(t, \mathbf{p}_0) = 0$, then Algorithm~\ref{alg:SA} terminates instantly. If $\delta(t, \mathbf{p}_0) > 0$, then Lemma~\ref{lem:rho_lower_incr} gives $\rho_{\mathrm{lower}}(t, \mathbf{p}_1) > \rho_{\mathrm{lower}}(t, \mathbf{p}_0) \geq 0$. Carrying on our reasoning in this way, we conclude that either Algorithm~\ref{alg:SA} terminates in a finite number of iteration or
       \begin{equation*}
           \rho_{\mathrm{lower}}(t, \mathbf{p}_{m + 1}) > \rho_{\mathrm{lower}}(t, \mathbf{p}_{m}) > 0.
       \end{equation*}
       for any $m \in \mathbb{N}$. According to Lemma~\ref{lem:upper_lower}, the sequence $\{\rho_{\mathrm{lower}}(t, \mathbf{p}_{m})\}_{m = 0}^{\infty}$ is upper bounded by $\rho_{\mathrm{upper}}(t, \mathbf{p})$. Moreover,
       \begin{equation*}
           0 < \Delta \overset{\mathrm{def}}{=} \alpha\rho_{\mathrm{lower}}(t, \mathbf{p}_1) \leq \alpha\rho_{\mathrm{lower}}(t, \mathbf{p}_m).
       \end{equation*}
       
       Thus, $\{\rho_{\mathrm{lower}}(t, \mathbf{p}_{m})\}_{m = 0}^{\infty}$ is a convergent sequence. Using
       \begin{equation*}
           \mathbf{p}_{m}(\boldsymbol{s}_{\mathcal{R}(t)}(\mathbf{p}_{m}) - \boldsymbol{s}_{\mathcal{R}(t)}(\mathbf{p}_{m + 1})) \geq 0,
       \end{equation*}
       \begin{equation*}
           (-\mathbf{p}_{m})(\boldsymbol{s}_{\mathcal{G}(t)}(-\mathbf{p}_{m}) - \boldsymbol{s}_{\mathcal{G}(t)}(-\mathbf{p}_{m + 1})) \geq 0,
       \end{equation*}
       \begin{multline*}
           \xi'(\gamma/2^{N(t, \mathbf{p}_m) - 1}; t, \mathbf{p}_m) \\
           = \boldsymbol{q}(t, \mathbf{p}_m)(\boldsymbol{s}_{\mathcal{G}(t)}(-\mathbf{p}_{m + 1}) - \boldsymbol{s}_{\mathcal{R}(t)}(\mathbf{p}_{m + 1})) \geq 0,
       \end{multline*}
       and
       \begin{equation*}
           \lVert\mathbf{p}_m + \gamma/2^{N(t, \mathbf{p}_m) - 1}\boldsymbol{q}(t, \mathbf{p}_m)\rVert < 1,
       \end{equation*}
       we obtain
       \begin{multline*}
           \rho_{\mathrm{lower}}(t, \mathbf{p}_{m + 1}) - \rho_{\mathrm{lower}}(t, \mathbf{p}_m)\\
           = \mathbf{p}_{m + 1}(\boldsymbol{s}_{\mathcal{G}(t)}(-\mathbf{p}_{m + 1}) - \boldsymbol{s}_{\mathcal{R}(t)}(\mathbf{p}_{m + 1})) \\
           - \mathbf{p}_{m}(\boldsymbol{s}_{\mathcal{G}(t)}(-\mathbf{p}_m) - \boldsymbol{s}_{\mathcal{R}(t)}(\mathbf{p}_{m}))\\
           = (\mathbf{p}_{m + 1} - \mathbf{p}_{m})(\boldsymbol{s}_{\mathcal{G}(t)}(-\mathbf{p}_{m + 1}) - \boldsymbol{s}_{\mathcal{R}(t)}(\mathbf{p}_{m + 1}))\\
           + \mathbf{p}_{m}(\boldsymbol{s}_{\mathcal{R}(t)}(\mathbf{p}_{m}) - \boldsymbol{s}_{\mathcal{R}(t)}(\mathbf{p}_{m + 1}))\\
           + (-\mathbf{p}_{m})(\boldsymbol{s}_{\mathcal{G}(t)}(-\mathbf{p}_{m}) - \boldsymbol{s}_{\mathcal{G}(t)}(-\mathbf{p}_{m + 1}))\\
           \geq (\mathbf{p}_{m + 1} - \mathbf{p}_{m})(\boldsymbol{s}_{\mathcal{G}(t)}(-\mathbf{p}_{m + 1}) - \boldsymbol{s}_{\mathcal{R}(t)}(\mathbf{p}_{m + 1}))\\
           = (\mathbf{p}_{m + 1} - (\mathbf{p}_m + \gamma/2^{N(t, \mathbf{p}_m) - 1}\boldsymbol{q}(t, \mathbf{p}_m)))\\
           \cdot(\boldsymbol{s}_{\mathcal{G}(t)}(-\mathbf{p}_{m + 1}) - \boldsymbol{s}_{\mathcal{R}(t)}(\mathbf{p}_{m + 1}))\\
           + \gamma/2^{N(t, \mathbf{p}_m) - 1}\boldsymbol{q}(t, \mathbf{p}_m)(\boldsymbol{s}_{\mathcal{G}(t)}(-\mathbf{p}_{m + 1}) - \boldsymbol{s}_{\mathcal{R}(t)}(\mathbf{p}_{m + 1}))\\
           \geq (\mathbf{p}_{m + 1} - (\mathbf{p}_m + \gamma/2^{N(t, \mathbf{p}_m) - 1}\boldsymbol{q}(t, \mathbf{p}_m)))\\
           \cdot(\boldsymbol{s}_{\mathcal{G}(t)}(-\mathbf{p}_{m + 1}) - \boldsymbol{s}_{\mathcal{R}(t)}(\mathbf{p}_{m + 1}))\\
           = (1 - \lVert\mathbf{p}_m + \gamma/2^{N(t, \mathbf{p}_m) - 1}\boldsymbol{q}(t, \mathbf{p}_m)\rVert)\rho_{\mathrm{lower}}(t, \mathbf{p}_{m + 1})\\
           > (1 - \lVert\mathbf{p}_m + \gamma/2^{N(t, \mathbf{p}_m) - 1}\boldsymbol{q}(t, \mathbf{p}_m)\rVert)\rho_{\mathrm{lower}}(t, \mathbf{p}_1) > 0.
       \end{multline*}
       Therefore,
       \begin{equation*}
           \lim_{m \to \infty}\lVert\mathbf{p}_m + \gamma/2^{N(t, \mathbf{p}_m) - 1}\boldsymbol{q}(t, \mathbf{p}_m)\rVert = 1.
       \end{equation*}
       Using~\eqref{eq:til_p_norm}, we conclude that
       \begin{equation*}
           \lim_{m \to \infty}\gamma/2^{N(t, \mathbf{p}_m) - 1}(1 - \gamma/2^{N(t, \mathbf{p}_m) - 1})\frac{\delta(t, \mathbf{p}_m)}{\rho_{\mathrm{upper}}(t, \mathbf{p}_m)} = 0.
       \end{equation*}
       If $\delta(t, \mathbf{p}_m) \leq \Delta$, then Algorithm~\ref{alg:SA} terminates. Taking into account
       \begin{equation*}
           \rho_{\mathrm{upper}}(t, \mathbf{p}_m) \leq \max_{\mathbf{s} \in \mathcal{R}(t), \: \tilde{\mathbf{s}} \in \mathcal{G}(t)}\lVert\tilde{\mathbf{s}} - \mathbf{s}\rVert < +\infty,
       \end{equation*}
       $\gamma/2^{N(t, \mathbf{p}_m) - 1} \leq \gamma < 1$, and $\delta(t, \mathbf{p}_m) > \Delta$, we obtain
       \begin{equation*}
           \lim_{m \to \infty}\gamma/2^{N(t, \mathbf{p}_m)} = 0, \quad  \lim_{m \to \infty}N(t, \mathbf{p}_m) = \infty.
       \end{equation*}
       According to definition~\eqref{eq:N_def}, we have
       \begin{equation*}
           \xi'(\gamma/2^{N(t, \mathbf{p}_m) - 1}; t, \mathbf{p}_m) \geq 0,
       \end{equation*}
       but $\xi'(\gamma/2^{N(t, \mathbf{p}_m) - 2}; t, \mathbf{p}_m) < 0$. Thus,
       \begin{multline*}
           0 > \xi'(\gamma/2^{N(t, \mathbf{p}_m) - 2}; t, \mathbf{p}_m) \\
           = \boldsymbol{q}(t, \mathbf{p}_m)(\boldsymbol{s}_{\mathcal{G}(t)}(-(\mathbf{p}_m + \gamma/2^{N(t, \mathbf{p}_m) - 2}\boldsymbol{q}(t, \mathbf{p}_m)))\\
           - \boldsymbol{s}_{\mathcal{R}(t)}(\mathbf{p}_m + \gamma/2^{N(t, \mathbf{p}_m) - 2}\boldsymbol{q}(t, \mathbf{p}_m)))\\
           = \boldsymbol{q}(t, \mathbf{p}_m)(\boldsymbol{s}_{\mathcal{G}(t)}(-(\mathbf{p}_m + \gamma/2^{N(t, \mathbf{p}_m) - 2}\boldsymbol{q}(t, \mathbf{p}_m)))\\
           - \boldsymbol{s}_{\mathcal{G}(t)}(-\mathbf{p}_m)) + \boldsymbol{q}(t, \mathbf{p}_m)(\boldsymbol{s}_{\mathcal{R}(t)}(\mathbf{p}_m)\\
           - \boldsymbol{s}_{\mathcal{R}(t)}(\mathbf{p}_m + \gamma/2^{N(t, \mathbf{p}_m) - 2}\boldsymbol{q}(t, \mathbf{p}_m))) + \delta(t, \mathbf{p}_m)\\
       \end{multline*}
       Calculating the limit and taking into account Lemma~\ref{lem:cont_cont}, we obtain
       \begin{equation*}
           \lim_{m \to \infty}\delta(t, \mathbf{p}_m) \leq 0.
       \end{equation*}
       It contradicts to $\delta(t, \mathbf{p}_m) > \Delta > 0$. Thus, there exists $m \in \mathbb{N}$ such that $\delta(t, \mathbf{p}_m) \leq \Delta \leq \alpha\rho_{\mathrm{lower}}(t, \mathbf{p}_m)$. It means that Algorithm~\ref{alg:SA} terminates in a finite number of iterations. Moreover,
       \begin{equation*}
           \lim_{\alpha \to +0} \delta(t, \mathbf{SA}_\alpha(t, \mathbf{p}; \gamma)) = 0.
       \end{equation*}
       Lemma~\ref{lem:limit_root} of Appendix~\ref{ap:additional} completes the proof.
    \end{pf}

    \begin{cor}\label{cor:sa_sep}
        $\rho_{\mathrm{lower}}(t, \tilde{\mathbf{p}}) > 0$ for $\tilde{\mathbf{p}} = \mathbf{SA}_{\alpha}(t, \mathbf{p}; \gamma)$.
    \end{cor}
    \begin{pf}
        Algorithm~\ref{alg:SA} terminates with
        \begin{equation*}
            \delta(t, \tilde{\mathbf{p}}) \leq \alpha\rho_{\mathrm{lower}}(t, \tilde{\mathbf{p}}).
        \end{equation*}
        Using~$\mathcal{R}(t) \cap \mathcal{G}(t) = \varnothing$, we obtain $\rho_{\mathrm{upper}}(t, \tilde{\mathbf{p}}) \geq \rho(t) > 0$. Thus,
        \begin{equation*}
            (1 + \alpha)\rho_{\mathrm{lower}}(t, \tilde{\mathbf{p}}) \geq \rho_{\mathrm{upper}}(t, \tilde{\mathbf{p}}) > 0.
        \end{equation*}
    \end{pf}

    \subsection{Gradient ascent distance algorithm}

    We now take advantage of the fact that the gradient of $\rho_{\mathrm{lower}}(t, \cdot)$ is described by Lemma~\ref{lem:lower_diff_p}. Thus, we can use the gradient ascent method to maximize $\rho_{\mathrm{lower}}(t, \cdot)$.
    
    \begin{algorithm}
    	\caption{Gradient ascent algorithm: $\mathbf{GA}_{\alpha}(t, \mathbf{p}; \gamma)$}
    	\label{alg:GA}
    	\textbf{Input:} $\alpha \in \mathbb{R}^+$, $t \in \mathbb{R}^+_0$, $\mathbf{p} \in \mathcal{S}^*$, $\gamma \in (0, 1)$\\
        \textbf{Require:} $\mathrm{RP1}$, $\mathrm{RT1}$, $\lVert\mathbf{p}\rVert = 1$, $\rho_{\mathrm{lower}}(t, \mathbf{p}) > 0$
    	\begin{algorithmic}[1]
            \While{$\delta(t, \mathbf{p}) > \alpha\rho_{\mathrm{lower}}(t, \mathbf{p})$}
                \While{$\rho_{\mathrm{lower}}(t, \mathbf{p} + \gamma\frac{\partial\rho_{\mathrm{lower}}(t, \mathbf{p})}{\partial\mathbf{p}}) \leq \rho_{\mathrm{lower}}(t, \mathbf{p})$}
                    \State $\gamma \gets \gamma / 2$
                \EndWhile
                \State $\mathbf{p} \gets \frac{\mathbf{p} + \gamma\frac{\partial\rho_{\mathrm{lower}}(t, \mathbf{p})}{\partial\mathbf{p}}}{\lVert\mathbf{p} + \gamma\frac{\partial\rho_{\mathrm{lower}}(t, \mathbf{p})}{\partial\mathbf{p}}\rVert}$
            \EndWhile
    	\end{algorithmic}	
    	\textbf{Output:} $\mathbf{p}$
    \end{algorithm}
    
    We will not prove the convergence of Algorithm~\ref{alg:GA}. Instead, we investigate the properties of the algorithm by numerical experiments in the following sections.

    \section{Algorithms of computing the optimal control}\label{sec:MTPLS}

    Now we describe algorithms for estimating parameters $T^*$, $\mathbf{p}^*$ of optimal control $\boldsymbol{u}_E(\boldsymbol{p}(\cdot; T^*, \mathbf{p}^*))$. Each algorithm uses GJK distance algorithm (Algorithm~\ref{alg:GJK}) for calculating the initial value of the support vector that corresponds to a strictly separative hyperplane.

    \subsection{Neustadt-Eaton algorithm}

    Algorithm~\ref{alg:NEB} draws upon ideas from the works of  \citet{Neustadt1960-ki} and \citet{Eaton1962-ai}. According to the general approach of \citet{Neustadt1960-ki}, the MTPLS is reduced to the maximization of the boosting-time function $F$. \citet{Eaton1962-ai} proposed to conduct the maximization in two stages: the inclination of the tangent hyperplane by tuning $\mathbf{p}$ (see Fig.~\ref{fig:incline}) and the time increment by the boosting-time function (see Fig.~\ref{fig:F}). 

    Algorithm~\ref{alg:NEB} is a steepest ascent algorithm where $\boldsymbol{q}(t, \mathbf{p})$ is the ascent direction for $F(t, \cdot)$ at $\mathbf{p}$.

    Line 6 of Algorithm~\ref{alg:NEB} is a constructive rule for choosing a step size that was proposed by~\citet{Eaton1962-ai}. In the general case, Algorithm~\ref{alg:NEB} can loop as there exists no proof of its convergence. \citet{Boltyanskii1971-ty} provided constructively verifiable conditions for checking the step size to guarantee the convergence:
    \begin{equation*}
        \rho_{\mathrm{lower}}(t, \mathbf{p} + \gamma\boldsymbol{q}(t, \mathbf{p})) \leq \gamma\rho_{\mathrm{upper}}^2(t, \mathbf{p}).
    \end{equation*}
    This inequality must replace the inequality from Line 6 in Algorithm~\ref{alg:NEB} to be able to prove its convergence. Unfortunately, it is not efficient to use this condition in practice, because it leads to a large number of iterations\footnote{We verified this using computational experiments with the isotropic rocket benchmark. An intuitive explanation for this phenomenon is as follows. After time increment with the function $F$, the lower estimate for the distance is zero when the upper estimate is non-zero. Taking into account the continuity of the lower estimate, we should expect its small changes at small step sizes. Therefore, the step size requires a sufficiently large reduction such that the product of the step length by the upper estimate becomes smaller than the corresponding lower estimate.}.
    
    \begin{algorithm}
    	\caption{Neustadt-Eaton algorithm: $\mathbf{NE}_{\varepsilon}(\mathbf{p})$}
    	\label{alg:NEB}
    	\textbf{Input:} $\varepsilon \in \mathbb{R}^+$, $\mathbf{p} \in \mathcal{S}^*$\\
        \textbf{Require:} RP1, RT2, $\lVert\mathbf{p}\rVert = 1$
    	\begin{algorithmic}[1]
            \State $t \gets 0$, $\gamma \gets 1$
            \If{$\lVert\mathbf{GJK}_{\varepsilon}(t, \mathbf{p})\rVert > \varepsilon$}
                \State $\mathbf{p} \gets -\mathbf{GJK}_{\varepsilon}(t, \mathbf{p})^\top/\lVert\mathbf{GJK}_{\varepsilon}(t, \mathbf{p})\rVert$
                \While{$\rho_{\mathrm{upper}}(t, \mathbf{p}) > \varepsilon$}
                    \State $t, \mathbf{p} \gets F(t, \mathbf{p}), \boldsymbol{p}(F(t, \mathbf{p}); t, \mathbf{p})$
                    \While{$\rho_{\mathrm{lower}}(t, \mathbf{p} + \gamma\boldsymbol{q}(t, \mathbf{p})) \leq \rho_{\mathrm{lower}}(t, \mathbf{p})$}
                        \State $\gamma \gets \gamma / 2$
                    \EndWhile
                    \State $\mathbf{p} \gets \frac{\mathbf{p} + \gamma\boldsymbol{q}(t, \mathbf{p})}{\lVert\mathbf{p} + \gamma\boldsymbol{q}(t, \mathbf{p})\rVert}$
                \EndWhile
            \EndIf
    	\end{algorithmic}	
    	\textbf{Output:} $t$, $\mathbf{p}$
    \end{algorithm}

    \subsection{Barr-Gilbert algorithm}

    Algorithm~\ref{alg:BG} is derived from ideas presented in the papers of \cite{Ho1962-ki}, \cite{Barr1969-aj}. According to the general approach of \cite{Ho1962-ki}, the MTPLS is decomposed to a sequence of the MDPs. We suppose that each MDP can be solved by some algorithm $\mathbf{DA}$. $\mathbf{DA}$ denotes the one of algorithms $\mathbf{GJK}^*$, $\mathbf{G}$, $\mathbf{SA}$, or $\mathbf{GA}$ for estimating the distance. So, $\mathbf{DA} \in \{\mathbf{GJK}^*, \mathbf{G}, \mathbf{SA}, \mathbf{GA}\}$. \cite{Barr1969-aj} proposed to use the boosting-time function $F(t, \mathbf{p})$ with $\mathbf{p}$ given by $\mathbf{DA}$ to construct an increasing sequence of lower estimations that converges to the optimal time.

    \begin{algorithm}
    	\caption{Barr-Gilbert algorithm: $\mathbf{BG}_{\varepsilon}(\mathbf{p}; \alpha)$}
    	\label{alg:BG}
    	\textbf{Input:} $\varepsilon \in \mathbb{R}^+$, $\mathbf{p} \in \mathcal{S}^*$, $\alpha \in \mathbb{R}^+$\\
        \textbf{Require:} RP1, RT2, $\lVert\mathbf{p}\rVert = 1$
    	\begin{algorithmic}[1]
            \State $t \gets 0$, $\gamma \gets 1$
            \If{$\lVert\mathbf{GJK}_{\varepsilon}(t, \mathbf{p})\rVert > \varepsilon$}
                \State $\mathbf{p} \gets -\mathbf{GJK}_{\varepsilon}(t, \mathbf{p})^\top/\lVert\mathbf{GJK}_{\varepsilon}(t, \mathbf{p})\rVert$
                \While{$\rho_{\mathrm{upper}}(t, \mathbf{p}) > \varepsilon$}
                    \State $t, \mathbf{p} \gets F(t, \mathbf{p}), \boldsymbol{p}(F(t, \mathbf{p}); t, \mathbf{p})$
                    \State $\mathbf{p} \gets \mathbf{DA}_{\alpha}(t, \mathbf{p}; \gamma)$
                \EndWhile
            \EndIf
    	\end{algorithmic}	
    	\textbf{Output:} $t$, $\mathbf{p}$
    \end{algorithm}

    \subsection{Semi-analytical algorithm}
    
    Algorithms~\ref{alg:NEB},~\ref{alg:BG} use the calculation of $F(t, \mathbf{p})$. Generally speaking, such calculations are associated with the need for numerical integration of the equations of dynamics until the function $\rho_{\mathrm{lower}}(\cdot, \boldsymbol{p}(\cdot, t, \mathbf{p}))$ changes sign. This procedure can be executed with explicit Runge-Kutta methods, for example. Note, that this approach is hard to adopt when attempting to use the analytical description of $\boldsymbol{s}_{\mathcal{R}(t)}$ to compute $F$.

    We now describe a novel algorithm that uses only the analytical description of  $\boldsymbol{s}_{\mathcal{R}(t)}$ and does not require the calculation of $F$. Algorithm~\ref{alg:S} alternates the procedure of time boosting using the function $f$ and the procedure of estimating the distance between the reachable set and the target set. We denote the output of Algorithm~\ref{alg:S} by $T^*_{\mathbf{S}}[\varepsilon, \alpha, \mathbf{p}]$, $\boldsymbol{p}^*_{\mathbf{S}}[\varepsilon, \alpha, \mathbf{p}]$.
    
    \begin{algorithm}
    	\caption{Semi-analytical algorithm: $\mathbf{S}_{\varepsilon}(\mathbf{p}; \alpha)$}
    	\label{alg:S}
    	\textbf{Input:} $\varepsilon \in \mathbb{R}^+$, $\mathbf{p} \in \mathcal{S}^*$, $\alpha \in \mathbb{R}^+$\\
        \textbf{Require:} $\mathrm{RP1}$, $\mathrm{RP2}$, $\mathrm{RT4}$, $\lVert\mathbf{p}\rVert = 1$
    	\begin{algorithmic}[1]
            \State $t \gets 0$, $\gamma \gets 1$
            \If{$\lVert\mathbf{GJK}_{\varepsilon}(t, \mathbf{p})\rVert > \varepsilon$}
                \State $\mathbf{p} \gets -\mathbf{GJK}_{\varepsilon}(t, \mathbf{p})^\top/\lVert\mathbf{GJK}_{\varepsilon}(t, \mathbf{p})\rVert$
                \While{$\rho_{\mathrm{upper}}(t, \mathbf{p}) > \varepsilon$}
                    \State $t \gets f(t, \mathbf{p})$
                    \State $\mathbf{p} \gets \mathbf{DA}_{\alpha}(t, \mathbf{p}; \gamma)$
                \EndWhile
            \EndIf
    	\end{algorithmic}	
    	\textbf{Output:} $t$, $\mathbf{p}$
    \end{algorithm}

    \begin{thm}
        Let $\mathrm{RP1}$, $\mathrm{RP2}$, $\mathrm{RT4}$, $\varepsilon \in \mathbb{R}^+$, $\alpha \in \mathbb{R}^+$, $\mathbf{p} \in \mathcal{S}^*$, $\lVert\mathbf{p}\rVert = 1$, $\mathbf{DA} \in \{\mathbf{GJK}^*, \mathbf{G}, \mathbf{SA}\}$\footnote{$\mathbf{GA}$ is excluded from this list since the distance algorithm has no proof of convergence}. If $T^* < +\infty$, then Algorithm~\ref{alg:S} terminates in a finite number of iterations and
        \begin{align*}
            &\mathbf{p}^* = \lim_{\varepsilon \to +0} \boldsymbol{p}^*_{\mathbf{S}}[\varepsilon, \alpha, \mathbf{p}],\\
            &T^* = \lim_{\varepsilon \to +0} T^*_{\mathbf{S}}[\varepsilon, \alpha, \mathbf{p}].
        \end{align*}
    \end{thm}
    \begin{pf}
        If $\lVert\mathbf{GJK}_{\varepsilon}(0, \mathbf{p})\rVert \leq \varepsilon$, then the distance between the plant and the target set is less than $\varepsilon$ at the initial time moment and Algorithm~\ref{alg:S} terminates instantly. Further, we suppose $\lVert\mathbf{GJK}_{\varepsilon}(0, \mathbf{p})\rVert > \varepsilon$. Let $\gamma_1$, $\gamma_2$, ... be a sequence of values of variable $\gamma$ in the iterations of while-loop. Let
        \begin{align*}
            &\mathbf{p}_1 \overset{\mathrm{def}}{=} -\frac{\mathbf{GJK}_{\varepsilon}(0, \mathbf{p})^\top}{\lVert\mathbf{GJK}_{\varepsilon}(0, \mathbf{p})\rVert},\quad \mathbf{p}_{n + 1} \overset{\mathrm{def}}{=} \mathbf{DA}_{\alpha}(t_n, \mathbf{p}_n; \gamma_n);\\
            &t_0 \overset{\mathrm{def}}{=} 0, \quad t_{n + 1} \overset{\mathrm{def}}{=} f(t_n, \mathbf{p}_{n + 1}).
        \end{align*}
        Using Corollary~\ref{cor:gjk_sep}, we deduce $\rho_{\mathrm{lower}}(t_0, \mathbf{p}_1) > 0$. Lemma~\ref{lem:f_separ} gives $\rho_{\mathrm{lower}}(t_1, \mathbf{p}_1) \geq 0$ and $t_1 \leq T^*$. If $\mathcal{R}(t_1) \cap \mathcal{G}(t_1) = \varnothing$, then $t_1 < T^*$ and Corollaries \ref{cor:gjks_sep}, \ref{cor:g_sep}, \ref{cor:sa_sep} give $\rho_{\mathrm{lower}}(t_1, \mathbf{p}_2) > 0$. If $\mathcal{R}(t_1) \cap \mathcal{G}(t_1) \neq \varnothing$, then $\rho_{\mathrm{lower}}(t_1, \tilde{\mathbf{p}}) \leq 0$ for any $\tilde{\mathbf{p}} \in \mathcal{S}^*$, $\tilde{\mathbf{p}} \neq \boldsymbol{0}$. Hence, $\rho_{\mathrm{lower}}(t_1, \mathbf{p}_1) = 0$. It implies $\boldsymbol{s}_{\mathcal{R}(t_1)}(\mathbf{p}_1) = \boldsymbol{s}_{\mathcal{G}(t_1)}(-\mathbf{p}_1)$ and $\rho_{\mathrm{upper}}(t_1, \mathbf{p}_1) = 0$. In this case, $\mathbf{p}_2 = \mathbf{p}_1$ and Algorithm~\ref{alg:S} terminates, because $\rho_{\mathrm{upper}}(t_1, \mathbf{p}_2) = 0 \leq \varepsilon$. Moreover, $t_1 = T^*$ and $\mathbf{p}_1 = \mathbf{p}^*$. Carrying on our reasoning in this way, we conclude that either Algorithm~\ref{alg:S} terminates in a finite number of iterations or $\rho_{\mathrm{lower}}(t_n, \mathbf{p}_{n + 1}) > 0$, $\rho_{\mathrm{upper}}(t_n, \mathbf{p}_{n + 1}) > \varepsilon$, and $t_n < T^*$ for all $n \in \mathbb{N}$. Now we prove that the second option is unavailable if $T^* < +\infty$.

        In the second option, $\{t_n\}_{n = 0}^\infty$ is an increasing sequence, because
        \begin{equation*}
            t_{n + 1} = t_{n} + \frac{\rho_{\mathrm{lower}}(t_n, \mathbf{p}_{n + 1})}{v_{\mathcal{R}} + v_{\mathcal{G}}} > t_{n}.
        \end{equation*}
        Thus, the sequence is convergent and
        \begin{equation*}
            \lim_{n \to \infty}\rho_{\mathrm{lower}}(t_n, \mathbf{p}_{n + 1}) = 0.
        \end{equation*}
        On the other hand, $\mathbf{DA} \in \{\mathbf{GJK}^*, \mathbf{G}, \mathbf{SA}\}$, $\mathbf{p}_{n + 1} = \mathbf{DA}_{\alpha}(t_n, \mathbf{p}_n; \gamma_n)$ and
        \begin{equation*}
            \delta(t_n, \mathbf{p}_{n + 1}) \leq \alpha\rho_{\mathrm{lower}}(t_n, \mathbf{p}_{n + 1}).
        \end{equation*}
        Hence,
        \begin{equation*}
            \varepsilon < \rho_{\mathrm{upper}}(t_n, \mathbf{p}_{n + 1}) \leq (1 + \alpha)\rho_{\mathrm{lower}}(t_n, \mathbf{p}_{n + 1}).
        \end{equation*}
        Taking a limit, we observe a contradiction. Thus, Algorithm~\ref{alg:S} terminates in a finite number of iterations for $T^* < +\infty$.

        Algorithm~\ref{alg:S} produces $T^*_{\mathbf{S}}[\varepsilon, \alpha, \mathbf{p}]$, $\boldsymbol{p}^*_{\mathbf{S}}[\varepsilon, \alpha, \mathbf{p}]$ as output and
        \begin{equation*}
            \rho_{\mathrm{upper}}(T^*_{\mathbf{S}}[\varepsilon, \alpha, \mathbf{p}], \boldsymbol{p}^*_{\mathbf{S}}[\varepsilon, \alpha, \mathbf{p}]) \leq \varepsilon.
        \end{equation*}
        It is clear that
        \begin{equation*}
             T^*_{\mathbf{S}}[\varepsilon_2, \alpha, \mathbf{p}] \leq T^*_{\mathbf{S}}[\varepsilon_1, \alpha, \mathbf{p}] \leq T^*
        \end{equation*}
        for $\varepsilon_1, \varepsilon_2 \in \mathbb{R}^+$ and $\varepsilon_1 \leq \varepsilon_2$. Thus, the corresponding limit exists and we denote it by
        \begin{equation*}
            T_{\mathbf{S}} \overset{\mathrm{def}}{=} \lim_{\varepsilon \to +0} T^*_{\mathbf{S}}[\varepsilon, \alpha, \mathbf{p}].
        \end{equation*}
        Let's prove, that $T_{\mathbf{S}} = T^*$. Suppose a contrary, that $T_{\mathbf{S}} < T^*$. Choosing $\varepsilon \in \mathbb{R}^+$ sufficiently small, we obtain a contradiction with
        \begin{equation*}
            \varepsilon \geq \rho_{\mathrm{upper}}(T^*_{\mathbf{S}}[\varepsilon, \alpha, \mathbf{p}], \boldsymbol{p}^*_{\mathbf{S}}[\varepsilon, \alpha, \mathbf{p}]) \geq \min_{t \in [0, T_{\mathbf{S}}]}\rho(t) > 0.
        \end{equation*}
        Thus,
        \begin{equation*}
            T^* = \lim_{\varepsilon \to +0} T^*_{\mathbf{S}}[\varepsilon, \alpha, \mathbf{p}].
        \end{equation*}
        Lemma~\ref{lem:upper_lower} gives
        \begin{multline*}
            \varepsilon \geq \rho_{\mathrm{upper}}(T^*_{\mathbf{S}}[\varepsilon, \alpha, \mathbf{p}], \boldsymbol{p}^*_{\mathbf{S}}[\varepsilon, \alpha, \mathbf{p}])\\
            -\rho_{\mathrm{lower}}(T^*_{\mathbf{S}}[\varepsilon, \alpha, \mathbf{p}], \boldsymbol{p}^*_{\mathbf{S}}[\varepsilon, \alpha, \mathbf{p}])\geq 0.
        \end{multline*}
        Hence,
        \begin{equation*}
            \lim_{\varepsilon \to +0}\delta(T^*_{\mathbf{S}}[\varepsilon, \alpha, \mathbf{p}], \boldsymbol{p}^*_{\mathbf{S}}[\varepsilon, \alpha, \mathbf{p}]) = 0.
        \end{equation*}
        Lemma~\ref{lem:delta_cont} gives
        \begin{equation*}
            \lim_{\varepsilon \to +0}\delta(T^*, \boldsymbol{p}^*_{\mathbf{S}}[\varepsilon, \alpha, \mathbf{p}]) = 0.
        \end{equation*}
        Finally, Lemma~\ref{lem:limit_root} from Appendix~\ref{ap:additional} yields
        \begin{equation*}
            \mathbf{p}^* = \lim_{\varepsilon \to +0} \boldsymbol{p}^*_{\mathbf{S}}[\varepsilon, \alpha, \mathbf{p}].
        \end{equation*}
    \end{pf}

    \begin{rem}
        If $T^* = +\infty$, then Algorithm~\ref{alg:S} doesn't terminate and the variable $t$ is increasing without limit.
    \end{rem}

    \section{Numerical experiments}\label{sec:exp}

    We now study the performance of the algorithms presented above. First, we describe the computational problems that arise when computing $\boldsymbol{s}_{\mathcal{R}(t)}(\mathbf{p})$, $\boldsymbol{p}(t; T, \mathbf{p})$, $F(t, \mathbf{p})$. We will also describe some issues that arise due to imprecise arithmetic in computations. Then we will describe an example of MTPLS, which demonstrates the advantages of using the analytic computation capability of $\boldsymbol{s}_{\mathcal{R}(t)}(\mathbf{p})$. Next, we will compare the empirical complexities of all algorithms.

    \subsection{Computational aspects}

    Each of the presented algorithms contains functions that cannot be computed analytically in the general case. However, if it is possible to compute $\boldsymbol{s}_{\mathcal{R}(t)}(\mathbf{p})$, $\boldsymbol{p}(t; T, \mathbf{p})$, and $F(t, \mathbf{p})$, then $\rho_{\mathrm{lower}}(t, \mathbf{p})$, $\rho_{\mathrm{upper}}(t, \mathbf{p})$, $\delta(t, \mathbf{p})$, $\boldsymbol{q}(t, \mathbf{p})$, $f(t, \mathbf{p})$, $\xi'(\gamma; t, \mathbf{p})$ can be easily evaluated from them.

    All algorithms must compute $\boldsymbol{s}_{\mathcal{R}(t)}(\mathbf{p})$. If an analytical description is not available, we can first obtain $\boldsymbol{p}(0; t, \mathbf{p})$ by backward integrating of~\eqref{eq:adj} and then integrate the equations of dynamics~\eqref{eq:dyn} together with the conjugate system~\eqref{eq:adj} and the control input $\boldsymbol{u}_E$ expressed through the conjugate vector. We use Runge-Kutta method to implement the integration (see~Algorithm~\ref{alg:sR}).
    
    \begin{algorithm}
    	\caption{Approximation of $\boldsymbol{s}_{\mathcal{R}(t)}(\mathbf{p})$}
    	\label{alg:sR}
    	\textbf{Input:} $\tau \in \mathbb{R}^+$, $t \in \mathbb{R}^+_0$, $\mathbf{p} \in \mathcal{S}^*$
    	\begin{algorithmic}[1]
            \State $T \gets t$, $\mathbf{s} \gets \mathbf{s}_0$
            \While{$t > 0$}
                \State $h \gets \tau$ \textbf{if} $t \geq \tau$ \textbf{else} $t$
                \State $\mathbf{k}^*_1 \gets -\mathbf{p}\boldsymbol{A}(t)$
                \State $\mathbf{k}^*_2 \gets -\left(\mathbf{p} - \frac{h}2\mathbf{k}^*_1\right)\boldsymbol{A}\left(t - \frac{h}2\right)$
                \State $\mathbf{k}^*_3 \gets -\left(\mathbf{p} - \frac{h}2\mathbf{k}^*_2\right)\boldsymbol{A}\left(t - \frac{h}2\right)$
                \State $\mathbf{k}^*_4 \gets -\left(\mathbf{p} - h\mathbf{k}^*_3\right)\boldsymbol{A}\left(t - h\right)$
                \State $\mathbf{p} \gets \mathbf{p} - \frac{h}{6}(\mathbf{k}^*_1 + 2\mathbf{k}^*_2 + 2\mathbf{k}^*_3 + \mathbf{k}^*_4)$
                \State $t \gets t - h$
            \EndWhile
            \While{$t < T$}
                \State $h \gets \tau$ \textbf{if} $T - \tau \geq \tau$ \textbf{else} $T - t$
                \State $\mathbf{k}^*_1 \gets -\mathbf{p}\boldsymbol{A}(t)$
                \State $\mathbf{k}^*_2 \gets -\left(\mathbf{p} + \frac{h}2\mathbf{k}^*_1\right)\boldsymbol{A}\left(t + \frac{h}2\right)$
                \State $\mathbf{k}^*_3 \gets -\left(\mathbf{p} + \frac{h}2\mathbf{k}^*_2\right)\boldsymbol{A}\left(t + \frac{h}2\right)$
                \State $\mathbf{k}^*_4 \gets -\left(\mathbf{p} + h\mathbf{k}^*_3\right)\boldsymbol{A}\left(t + h\right)$
                \State $\mathbf{k}_1 \gets \boldsymbol{A}\left(t\right)\mathbf{s} + \boldsymbol{u}_E(\mathbf{p})$
                \State $\mathbf{k}_2 \gets \boldsymbol{A}\left(t + \frac{h}2\right)\left(\mathbf{s} + \frac{h}2\mathbf{k}_1\right) + \boldsymbol{u}_E(\mathbf{p} + \frac{h}2\mathbf{k}^*_1)$
                \State $\mathbf{k}_3 \gets \boldsymbol{A}\left(t + \frac{h}2\right)\left(\mathbf{s} + \frac{h}2\mathbf{k}_2\right) + \boldsymbol{u}_E(\mathbf{p} + \frac{h}2\mathbf{k}^*_2)$
                \State $\mathbf{k}_4 \gets \boldsymbol{A}(t + h)(\mathbf{s} + h\mathbf{k}_3) + \boldsymbol{u}_E(\mathbf{p} + h\mathbf{k}^*_3)$
                \State $\mathbf{p} \gets \mathbf{p} + \frac{h}{6}(\mathbf{k}^*_1 + 2\mathbf{k}^*_2 + 2\mathbf{k}^*_3 + \mathbf{k}^*_4)$
                \State $\mathbf{s} \gets \mathbf{s} + \frac{h}{6}(\mathbf{k}_1 + 2\mathbf{k}_2 + 2\mathbf{k}_3 + \mathbf{k}_4)$
                \State $t \gets t + h$
            \EndWhile
    	\end{algorithmic}	
    	\textbf{Output:} $\mathbf{s}$
    \end{algorithm}

    Algorithms~\ref{alg:NEB},~\ref{alg:BG} additionally require the calculation of $F(t, \mathbf{p})$, $\boldsymbol{p}(F(t, \mathbf{p}); t, \mathbf{p})$. In fact,
    \begin{equation}\label{eq:FpsR}
        F(t, \mathbf{p}), \: \boldsymbol{p}(F(t, \mathbf{p}); t, \mathbf{p}), \: \boldsymbol{s}_{\mathcal{R}(F(t, \mathbf{p}))}(\boldsymbol{p}(F(t, \mathbf{p}); t, \mathbf{p}))
    \end{equation}
    can be calculated simultaneously within one numerical integration, if $t$, $\mathbf{p}$, $\boldsymbol{s}_{\mathcal{R}(t)}(\mathbf{p})$ are known initially (see Algorithm~\ref{alg:FpsR}).
    
    \begin{algorithm}
    	\caption{Approximation of~\eqref{eq:FpsR}}
    	\label{alg:FpsR}
    	\textbf{Input:} $\tau \in \mathbb{R}^+$, $t \in \mathbb{R}^+_0$, $\mathbf{p} \in \mathcal{S}^*$, $\mathbf{s} \in \mathcal{S}$\\
        \textbf{Require:} $\mathbf{s} = \boldsymbol{s}_{\mathcal{R}(t)}(\mathbf{p})$, $\mathbf{p}(\boldsymbol{s}_{\mathcal{G}(t)}(-\mathbf{p}) - \mathbf{s}) > 0$
    	\begin{algorithmic}[1]
            \State $\tilde{t} \gets t$, $\tilde{\mathbf{p}} \gets \mathbf{p}$, $\tilde{\mathbf{s}} \gets \mathbf{s}$ 
            \While{$\mathbf{p}(\boldsymbol{s}_{\mathcal{G}(t)}(-\mathbf{p}) - \mathbf{s}) > 0$}
                \State $\tilde{t} \gets t$, $\tilde{\mathbf{p}} \gets \mathbf{p}$, $\tilde{\mathbf{s}} \gets \mathbf{s}$ 
                \State $\mathbf{k}^*_1 \gets -\mathbf{p}\boldsymbol{A}(t)$
                \State $\mathbf{k}^*_2 \gets -\left(\mathbf{p} + \frac{\tau}2\mathbf{k}^*_1\right)\boldsymbol{A}\left(t + \frac{\tau}2\right)$
                \State $\mathbf{k}^*_3 \gets -\left(\mathbf{p} + \frac{\tau}2\mathbf{k}^*_2\right)\boldsymbol{A}\left(t + \frac{\tau}2\right)$
                \State $\mathbf{k}^*_4 \gets -\left(\mathbf{p} + \tau\mathbf{k}^*_3\right)\boldsymbol{A}\left(t + \tau\right)$
                \State $\mathbf{k}_1 \gets \boldsymbol{A}\left(t\right)\mathbf{s} + \boldsymbol{u}_E(\mathbf{p})$
                \State $\mathbf{k}_2 \gets \boldsymbol{A}\left(t + \frac{\tau}2\right)\left(\mathbf{s} + \frac{\tau}2\mathbf{k}_1\right) + \boldsymbol{u}_E(\mathbf{p} + \frac{\tau}2\mathbf{k}^*_1)$
                \State $\mathbf{k}_3 \gets \boldsymbol{A}\left(t + \frac{\tau}2\right)\left(\mathbf{s} + \frac{\tau}2\mathbf{k}_2\right) + \boldsymbol{u}_E(\mathbf{p} + \frac{\tau}2\mathbf{k}^*_2)$
                \State $\mathbf{k}_4 \gets \boldsymbol{A}(t + \tau)(\mathbf{s} + \tau\mathbf{k}_3) + \boldsymbol{u}_E(\mathbf{p} + \tau\mathbf{k}^*_3)$
                \State $\mathbf{p} \gets \mathbf{p} + \frac{\tau}{6}(\mathbf{k}^*_1 + 2\mathbf{k}^*_2 + 2\mathbf{k}^*_3 + \mathbf{k}^*_4)$
                \State $\mathbf{s} \gets \mathbf{s} + \frac{\tau}{6}(\mathbf{k}_1 + 2\mathbf{k}_2 + 2\mathbf{k}_3 + \mathbf{k}_4)$
                \State $t \gets t + \tau$
            \EndWhile
    	\end{algorithmic}	
    	\textbf{Output:} $\tilde{t}$, $\tilde{\mathbf{p}}$, $\tilde{\mathbf{s}}$
    \end{algorithm}

    \subsection{Issues of imprecise arithmetic}
    
    Imprecise arithmetic appears due to numerical integrations or floating-point calculations. Imprecise arithmetic can lead the implementations of algorithms to infinite loop, or result in irrelevant solutions such that $\rho_{\mathrm{upper}}(t, \mathbf{p}) > \varepsilon$ for terminal $t$, $\mathbf{p}$. In the following, we investigate three algorithms $\mathbf{NE}$, $\mathbf{BG}$, $\mathbf{S}$. Algorithms $\mathbf{BG}$ and $\mathbf{S}$ use four distance algorithms: $\mathbf{GJK}^*$, $\mathbf{G}$, $\mathbf{SA}$, $\mathbf{GA}$.

    We use inequalities derived from theorems and lemmas about algorithms to detect issues with imprecise arithmetic. Below we provide a list of failure conditions:
    \begin{enumerate}
        \item $\delta(t, \mathbf{p}) < 0$;
        \item $\rho_{\mathrm{lower}}(t, \mathbf{p}) < 0$;
        \item $t$ is decreasing;
        \item $F$ is called more than $10^4$ times;
        \item $\lVert\mathbf{s}\rVert$ is non-decreasing;
        \item $\lVert\mathbf{s}\rVert = 0$;
        \item $\gamma \leq 0$.
    \end{enumerate}
    Conditions 1-4 are verified for all algorithms. Conditions 5, 6 are tested only for algorithms $\mathbf{BG}$, $\mathbf{S}$ that use $\mathbf{GJK}^*$ or $\mathbf{G}$ as distance algorithms. Condition 7 is tested for $\mathbf{NE}$, $\mathbf{BG}+\mathbf{SA}$, $\mathbf{BG}+\mathbf{GA}$, $\mathbf{S}+\mathbf{SA}$, $\mathbf{S}+\mathbf{GA}$.
    
    \subsection{Isotropic rocket benchmark}
    
    To examine the algorithms we use a model that describes the motion of a material point in a viscous medium under a force that is arbitrary in direction, but limited in magnitude. This model was named "isotropic rocket" in the early work of \cite{Isaacs1955-bn}.

    The state space for the model is $\mathcal{S} = \mathbb{R}^{4 \times 1}$. The dynamics of the isotropic rocket is described by~\eqref{eq:dyn} with
    \begin{equation*}
        \boldsymbol{A}(t) =
        \begin{bmatrix}
            0 & ~~0 & ~~1 & ~~0\\
            0 & ~~0 & ~~0 & ~~1\\
            0 & ~~0 & -1 & ~~0\\
            0 & ~~0 & ~~0 & -1
        \end{bmatrix}.
    \end{equation*}
    The initial state is $\mathbf{s}_0 = \begin{bmatrix}0 & 0 & v_0 & 0\end{bmatrix}^\top$, where $v_0 \in [0, 1)$ is a parameter of the model. The range of control inputs is
    \begin{equation*}
        \mathcal{U} = \left\{\mathbf{I}_{\frac12}\mathbf{u}:\: \mathbf{u} \in \mathcal{S}, \: \lVert\mathbf{I}_{\frac12}\mathbf{u}\rVert \leq 1\right\},
    \end{equation*}
    where
    \begin{equation*}
        \mathbf{I}_{\frac12} \overset{\mathrm{def}}{=}
        \begin{bmatrix}
            0 & 0 & 0 & 0\\
            0 & 0 & 0 & 0\\
            0 & 0 & 1 & 0\\
            0 & 0 & 0 & 1
        \end{bmatrix}.
    \end{equation*}
    
    The advantage of using the isotropic rocket as a benchmark is the ability to calculate analytically $\boldsymbol{p}(t; T, \mathbf{p})$ and $\boldsymbol{s}_{\mathcal{R}(t)}(\mathbf{p})$ \citep[Eqs.~2,~7,~8]{Buzikov2024-zm}. The isotropic rocket fulfills requirements $\mathrm{RP1}$, $\mathrm{RP2}$ and $v_{\mathcal{R}} = \sqrt{5}$. The maximum condition~\eqref{eq:max_princ} gives
    \begin{equation*}
        \boldsymbol{u}_E(\mathbf{p}) = \frac{\mathbf{I}_{\frac12}\mathbf{p}^\top}{\lVert\mathbf{p}\mathbf{I}_{\frac12}\rVert}.
    \end{equation*}
    We will use a point moving with constant velocity along a straight line as the target set:
    \begin{equation*}
        \mathcal{G}(t) \overset{\mathrm{def}}{=}
        \left\{
        \begin{bmatrix}
            s_1 + v_1 t & s_2 + v_2 t & v_1 & v_2
        \end{bmatrix}^\top\right\}.
    \end{equation*}
    Here, $s_1, s_2, v_1, v_2 \in \mathbb{R}$ are parameters of the target motion and $v_1^2 + v_2^2 < 1$. This motion fulfills requirement $\mathrm{RT4}$ with $v_{\mathcal{G}} = \sqrt{v_1^2 + v_2^2}$.
    
    We use the following ranges for parameters:
    \begin{align*}
        &\varepsilon \in \{3^{-i}:\: i \in \{3, 4, ..., 11\}\},\\
        &v_0 \in \left\{\frac{i}{8}:\: i \in \{0, 1, ..., 7\}\right\},\\
        &s_1 \in \left\{\frac{10i}{13} - 5:\: i \in \{0, 1, ..., 13\}\right\},\\
        &s_2 \in \left\{\frac{5i}{13}:\: i \in \{0, 1, ..., 13\}\right\},
    \end{align*}
    \begin{multline*}
        (v_1, v_2) \in \left\{\left(\frac{j}8\cos\frac{2\pi i}{10}, \frac{j}8\sin\frac{2\pi i}{10}\right):\: i \in \{0, 1, ..., 9\},\right.\\
        j \in \{0, 1, ..., 4\}\bigg\}.
    \end{multline*}
    Thus, the total number of samples is $9\cdot8\cdot14\cdot14\cdot10\cdot5=705600$. The fixed value of $\varepsilon$ corresponds to $78400$ samples.

    As an initial approximation for $\mathbf{p}$ in all cases, we will take the following:
    \begin{equation*}
        \mathbf{p} = \frac{\begin{bmatrix}
            s_1 & s_2 & v_1 - v_0 & v_2
        \end{bmatrix}}{\sqrt{s_1^2 + s_2^2 + (v_1 - v_0)^2 + v_2^2}}.
    \end{equation*}

    We use three values of numerical integration step-size $\tau \in \{10^{-2}, 10^{-3}, 10^{-4}\}$ in Algorithm~\ref{alg:FpsR}.

    \subsection{Algorithms' performance}
    
\begin{figure*}
        \centering
        \includegraphics[width=0.33\textwidth]{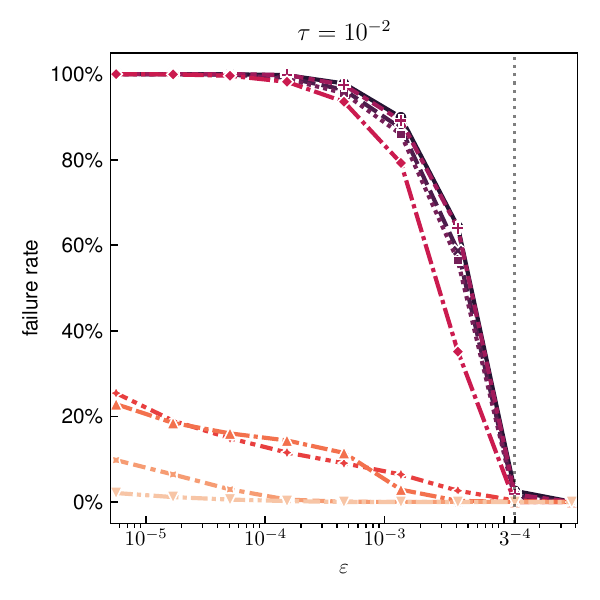}
        \includegraphics[width=0.33\textwidth]{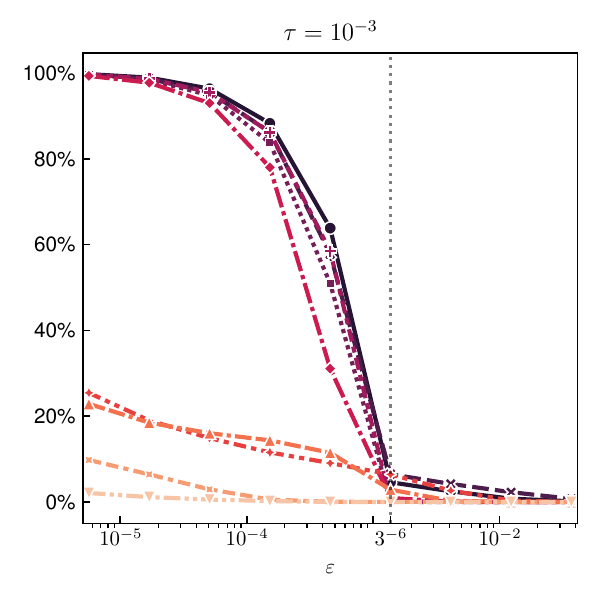}
        \includegraphics[width=0.33\textwidth]{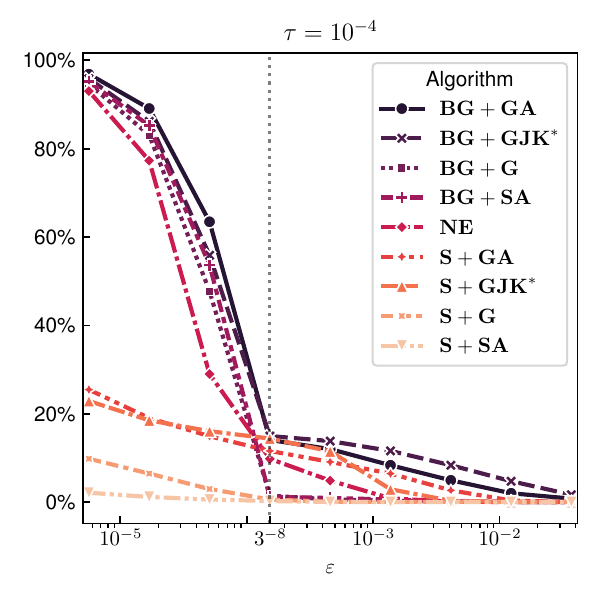}
        \caption{Dependence of the failure rate on the precision $\varepsilon$. The dotted vertical line is approximately located at the threshold value, after which there is a large growth of failure rate.}
        \label{fig:fail_rate}
    \end{figure*}

    First, we give statistics on the failure rate, i.e., the relative number of samples in which imprecise arithmetic yields failure conditions. Fig.~\ref{fig:fail_rate} illustrates a dramatic increase in the failure rate for algorithms that exploit the numerical integration ($\mathbf{NE}$, $\mathbf{BG}$ with all distance algorithms). The start of this growth occurs at $\varepsilon \approx \tau$ as expected. Note that algorithms that do not use the numerical integration ($\mathbf{S}$ with all distance algorithms) have a non-zero failure rate. This is explained by rounding errors for float64. $\mathbf{S} + \mathbf{SA}$ demonstrates the lowest failure rate.

    \begin{figure*}
        \centering
        \includegraphics[width=0.33\textwidth]{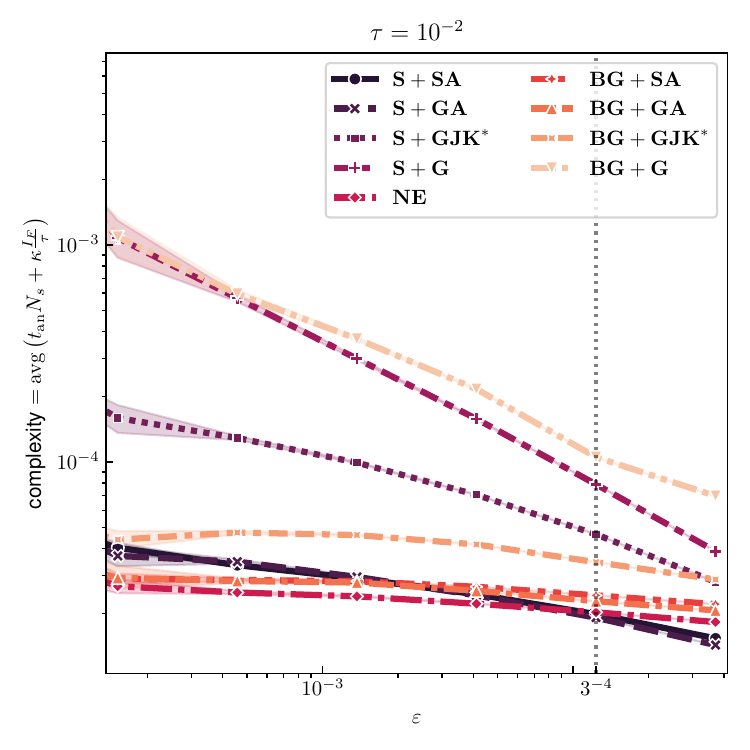}
        \includegraphics[width=0.33\textwidth]{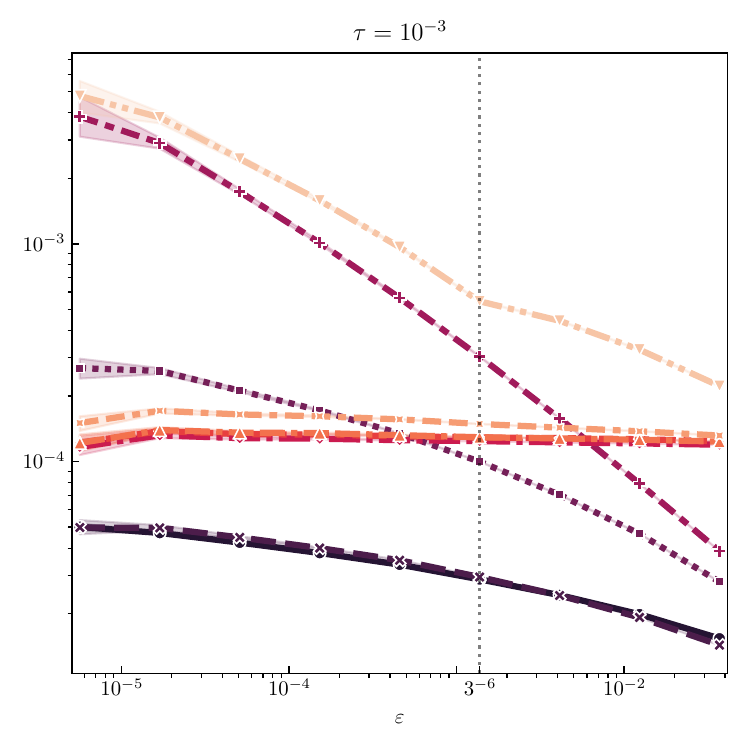}
        \includegraphics[width=0.33\textwidth]{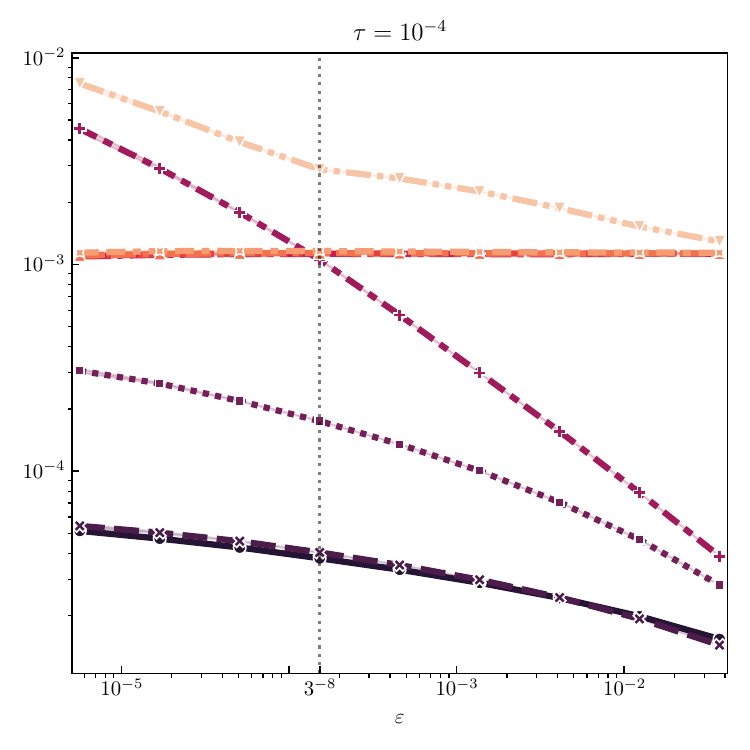}
        \caption{Dependence of Type B complexity on the precision $\varepsilon$. The dotted lines mark the thresholds of rapid increase in the failure rate (see~Fig.~\ref{fig:fail_rate}). If one of the algorithms had a failure in the sample, this sample was not taken into account in the averaging for all algorithms. Thus, the relevant data are to the right of the thresholds.}
        \label{fig:complexity_s}
    \end{figure*}

    We now turn to comparing the empirical complexities of the algorithms. We first note that the chosen design of the algorithms is such that the amount of memory allocated depends only on the dimension of the state space $\mathcal{S}$ and is independent of $\tau$. Therefore, all algorithms have the same memory complexity. This means that we should compare only the runtime complexity. Before defining the notion of complexity, we present a classification of MTPLS. Conventionally, we can categorize any MTPLS into one of three types:
    \begin{itemize}
        \item[Type A] $F$, $\boldsymbol{s}_{\mathcal{R}(\cdot)}$ are not available analytically and can be calculated only by numerical integrating.
        \item[Type B] $\boldsymbol{s}_{\mathcal{R}(\cdot)}$ is available analytically and $F$ can be calculated only by numerical integrating.
        \item[Type C] $F$, $\boldsymbol{s}_{\mathcal{R}(\cdot)}$ are available analytically.
    \end{itemize}
    The isotropic rocket benchmark is of Type B. But we are also interested in observing the complexity of the algorithms if we could suddenly forget that we can compute $\boldsymbol{s}_{\mathcal{R}(\cdot)}$ analytically (Type A), or, conversely, can compute $F$ analytically as well (Type C).

    We avoid the direct use of CPU-time as an empirical complexity measure because it is a subjective indicator that depends on the quality of the implementation. Instead, we use the following empirical complexity metrics:
    \begin{align*}
        &\text{Type A: complexity}= \mathrm{avg}\left(3I_s + 2I_F\right),\\
        &\text{Type B: complexity}= \mathrm{avg}\left(t_{\mathrm{an}}N_s + \kappa\frac{I_F}{\tau}\right),\\
        &\text{Type C: complexity}= \mathrm{avg}\left(N_s + N_F\right).
    \end{align*}
    The averaging is computed for each of $\varepsilon$ over all $v_0$, $s_1$, $s_2$, $v_1$, $v_2$. Here $I_s$, $I_F$ are the sums of the lengths of the intervals on which the numerical integration was performed by Algorithms~\ref{alg:sR},~\ref{alg:FpsR}; $\tau$ is a step-size of numerical integration; $N_s$, $N_F$ are the numbers of unique calls of $\boldsymbol{s}_{\mathcal{R}(\cdot)}$ and $F$; $t_{\mathrm{an}} \approx 422\cdot 10^{-9} \text{ sec.}$ is the average CPU-time of the call of $\boldsymbol{s}_{\mathcal{R}(\cdot)}$; $\kappa \approx 208 \cdot 10^{-10} \text{ sec.}$ is the average CPU-time for integrating with Algorithm~\ref{alg:FpsR} on a unit time interval $[0, 1]$ with $\tau = 1$. The values $t_{\mathrm{an}}$, $\kappa$ are measured on Intel Core i5-1135G7, Python 3.10.12, numba 0.60.0.

    Let us explain the meaning of choosing such formulas for empirical complexities. For Type A the main operation is numerical integrating by Algorithms~\ref{alg:sR},~\ref{alg:FpsR}. $I_s$ is associated with Algorithm~\ref{alg:sR} that requires a one backward integration and two forward integrations. $I_F$ is associated with Algorithm~\ref{alg:FpsR} that requires two forward integrations. That's why $I_s$ is taken with a multiplier of $3$ and $I_F$ is taken with a multiplier of $2$.
    
    For Type B\footnote{It would be fairer to calculate complexity by $\mathrm{avg}\left(\frac{t_{\mathrm{an}}N_s}{t} + \kappa\frac{I_F}{t\tau}\right)$, where $t$ is the optimal value of sample. This approach balances the averaging of both high and low optimal values of the samples. We have separately verified that such balancing has a negligible effect on the appearance of the plots.}, we should compare the cost of analytical calculation of $\boldsymbol{s}_{\mathcal{R}(\cdot)}$ with the cost of numerical integration to obtain $F$. The coefficients $t_{\mathrm{an}}$, $\kappa$ serve as weights to simultaneously account for the number of $\boldsymbol{s}_{\mathcal{R}(\cdot)}$-calls and the integration time $I_F$ required to compute the function $F$. Note that the empirical complexity for Type B is expressed in seconds, in contrast to the empirical complexity metrics for Types A, C. We can associate this amount of seconds with the CPU-time of the algorithm if the time was spent only on the computation of $\boldsymbol{s}_{\mathcal{R}(\cdot)}$, $F$.

    For Type C, we assume that analytically computing $\boldsymbol{s}_{\mathcal{R}(\cdot)}$ is as costly as computing $F$. Thus, the measure of empirical complexity is the total number of unique calls to these functions.

\begin{figure}
        \centering
        \includegraphics[width=0.43\textwidth]{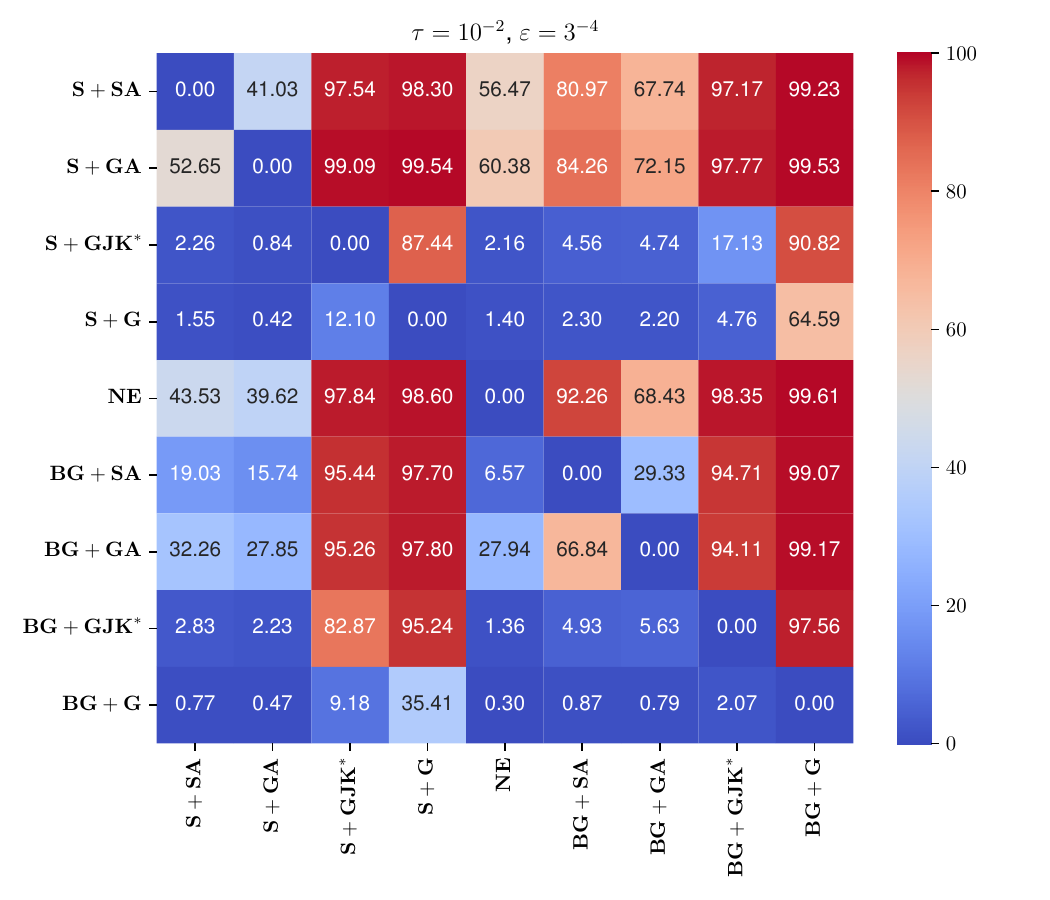}
        \includegraphics[width=0.43\textwidth]{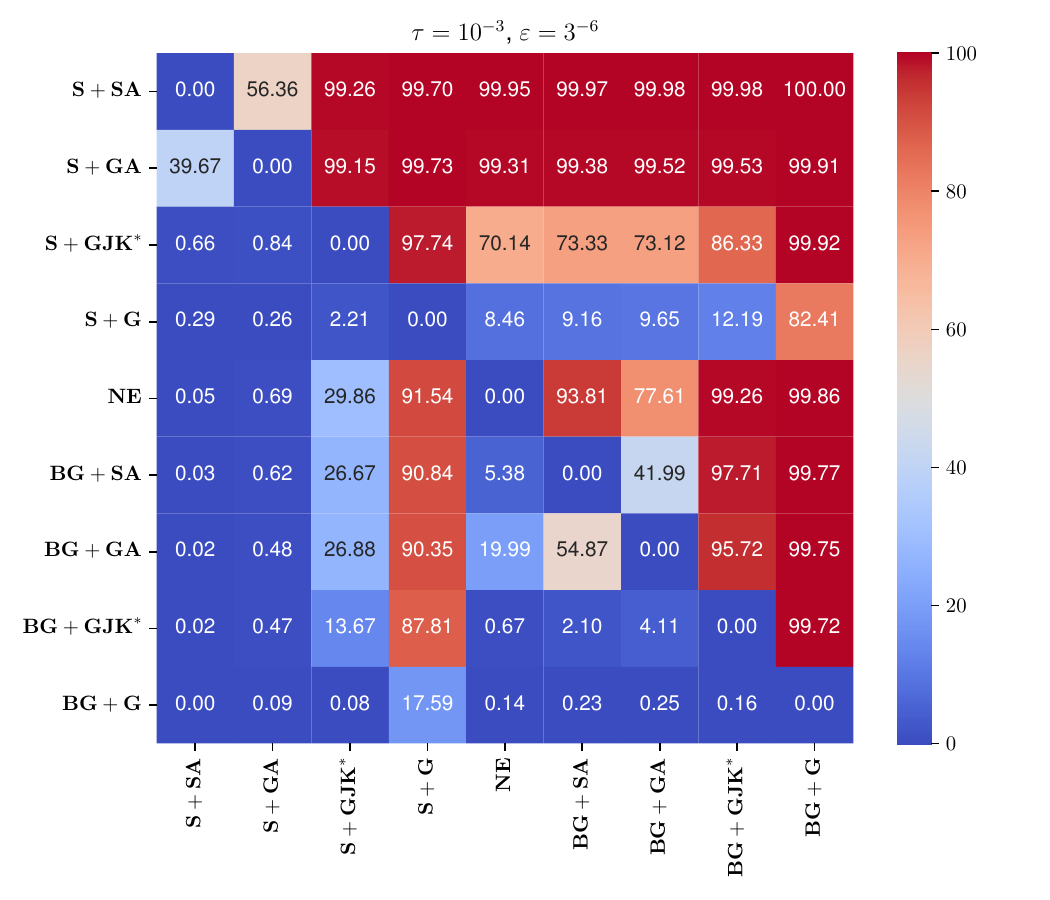}
        \includegraphics[width=0.43\textwidth]{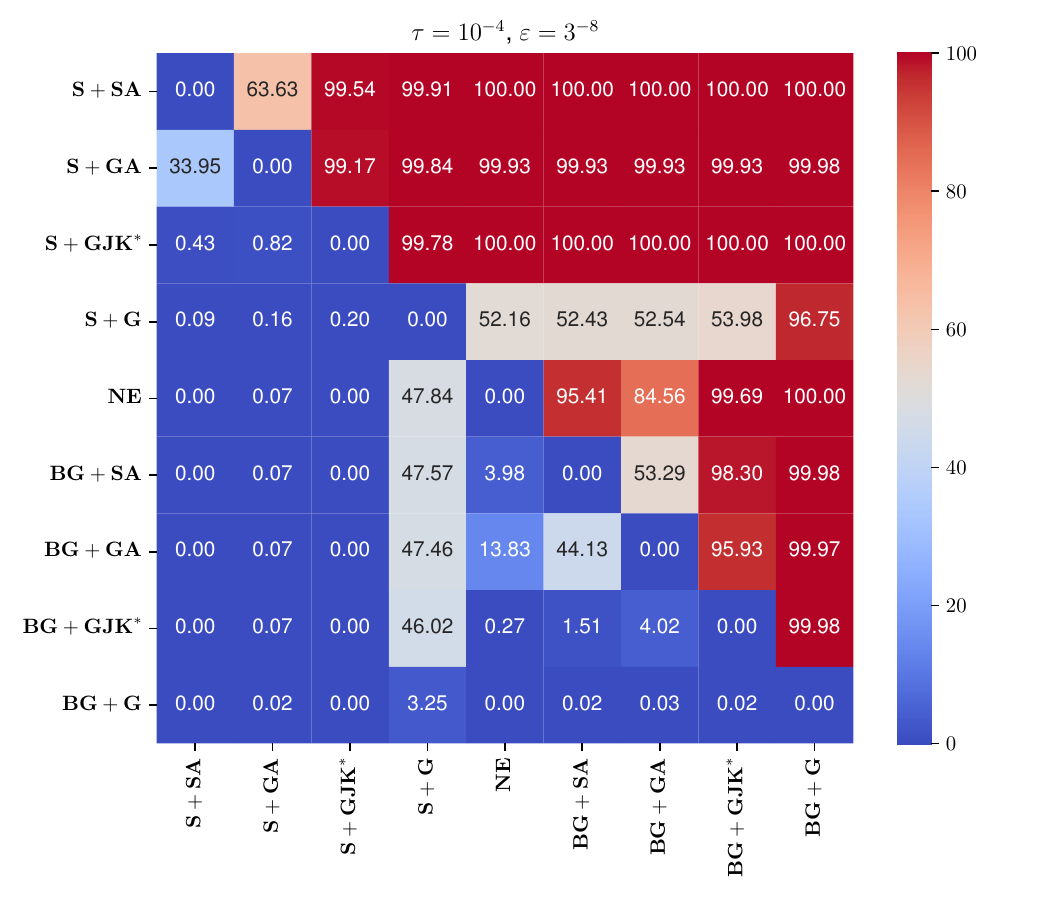}
        \caption{The superiority maps of algorithms for Type B complexity. Each cell contains the percentage of samples in which the row-naming algorithm demonstrates lower complexity than the column-naming algorithm. For example, $\mathbf{S}+\mathbf{SA}$ has lower complexity than $\mathbf{S} + \mathbf{GA}$ in $63.63\%$ of the samples for $\tau = 10^{-4}$, $\varepsilon=3^{-8}$. We excluded all samples with algorithm failures during computation.}
        \label{fig:superiority_map}
    \end{figure}

    A comparison of Type B empirical complexities is presented in Fig.~\ref{fig:complexity_s}. Fig.~\ref{fig:superiority_map} illustrates the superiority maps of the algorithms over other algorithms for different precisions and integration step-sizes. From the above data, we can conclude that $\mathbf{S}+\mathbf{SA}$ and $\mathbf{S} + \mathbf{GA}$ are superior to other algorithms. For high-precision computations, these algorithms outperform others by tens or hundreds of times (see Fig.~\ref{fig:complexity_s}).

\begin{figure}
        \centering
        \includegraphics[width=0.47\textwidth]{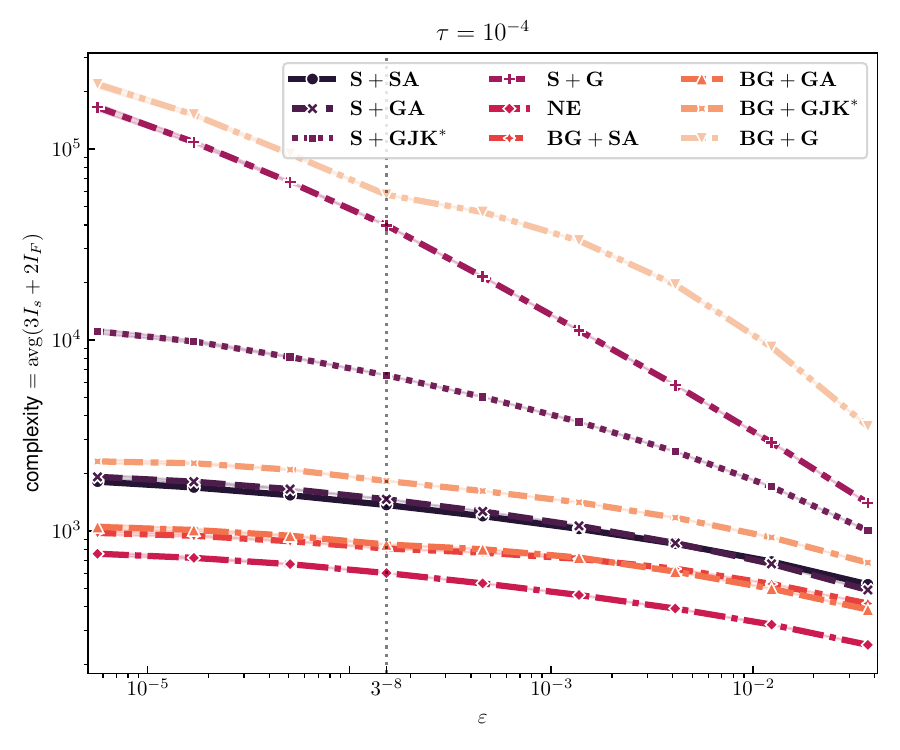}
        \caption{Dependence of Type A complexity on the precision $\varepsilon$}
        \label{fig:complexity_i}
    \end{figure}

\begin{figure}
        \centering
        \includegraphics[width=0.47\textwidth]{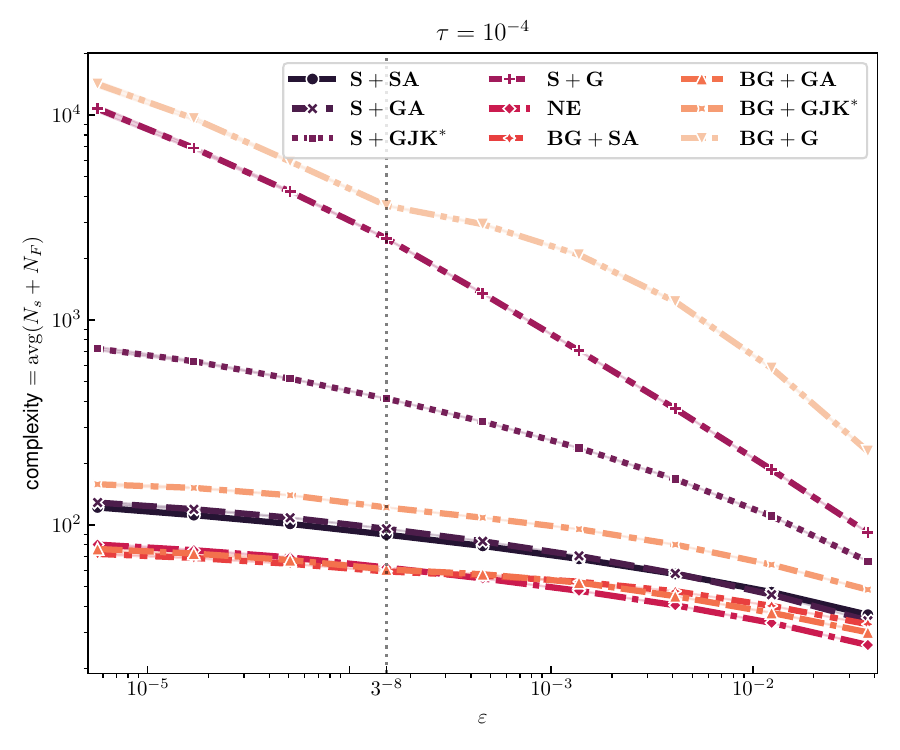}
        \caption{Dependence of Type C complexity on the precision $\varepsilon$}
        \label{fig:complexity_Fa}
    \end{figure}

    Figs.~\ref{fig:complexity_i},~\ref{fig:complexity_Fa} illustrate Type A and Type C complexities. From this data, we can conclude that $\mathbf{NE}$ is the most promising algorithm in Type A problem. Thus, if we were not able to compute $\boldsymbol{s}_{\mathcal{R}(\cdot)}$ analytically, we should have preferred the $\mathbf{NE}$ algorithm.
    
    We can also conclude that $\mathbf{NE}$, $\mathbf{BG}+\mathbf{SA}$, $\mathbf{BG}+\mathbf{GA}$ outperform other algorithms in Type C problem. Thus, if we were able to compute $F$ analytically as well as $\boldsymbol{s}_{\mathcal{R}(\cdot)}$, the preference should be given to $\mathbf{NE}$, $\mathbf{BG}+\mathbf{SA}$, $\mathbf{BG}+\mathbf{GA}$.

    \section{Conclusions}\label{sec:concl}

    This paper presents a contemporary and unified perspective on classical algorithms for solving minimum-time problems for linear systems. Furthermore, we introduced a novel algorithm and provided a constructive proof of its convergence. The results of the numerical experiments demonstrated that the novel algorithm exhibited a low failure rate and a relatively small empirical complexity in the context of the isotropic rocket benchmark. 
    
    It is noteworthy that the efficient solution of the problem for the isotropic rocket is of intrinsic importance. It is notable that there are few models in the existing literature that effectively address the problem of relocation to a given point with a given velocity vector. One such model is the Dubins model. We note that if the destination point is moving, then even in the case of the Dubins model, it cannot be said that an efficient solution to the problem is guaranteed in the general case. In the case of the isotropic rocket, we can now state that an effective solution to this problem has been achieved.

    It is anticipated that researchers interested in obtaining computationally efficient solutions to specific minimum-time problems for linear systems will direct their attention to the described algorithms. We would like to direct their attention to the Neustadt-Eaton algorithm as a tool that has demonstrated efficacy in situations where the reachable set cannot be parameterized analytically and all intermediate values must be obtained through numerical integration. Furthermore, it should be noted that when the geometry of the target set is complex, the Gilbert-Johnson-Keerthi algorithm is capable of identifying an admissible input for the Neustadt-Eaton algorithm.

    In the present study, we concentrated on the contemporary description of classical algorithms for the extended problem statement. We have supplied a constructive proof of convergence solely for the novel algorithm. The convergence of the Neustadt-Eaton and Barr-Gilbert algorithms has not been examined for the new extended problem statement. In our view, this investigation represents a promising continuation of this work.

    \begin{ack}
        The research was supported by RSF (project No. 23-19-00134).
    \end{ack}

    \appendix
    
    \section{Auxiliary statement}\label{ap:additional}

    \begin{lem}\label{lem:limit_root}
        Let $\mathcal{X}$ be a compact set in a finite dimension normed space, $f: \mathcal{X} \to \mathbb{R}$ be a continuous function, $\boldsymbol{\xi}: \mathbb{R} \to \mathcal{X}$ be an arbitrary function. Let $\mathbf{x}^* \in \mathcal{X}$ be a unique root of $f(\mathbf{x}) = 0$, i.e. $f(\mathbf{x}^*) = 0$ and $f(\mathbf{x}) \neq 0$ for any $\mathbf{x} \in \mathcal{X}$, $\mathbf{x} \neq \mathbf{x}^*$. If
        \begin{equation*}
            \lim_{\varepsilon \to 0} f(\boldsymbol{\xi}(\varepsilon)) = 0,
        \end{equation*}
        then
        \begin{equation*}
            \lim_{\varepsilon \to 0} \boldsymbol{\xi}(\varepsilon) = \mathbf{x}^*.
        \end{equation*}
    \end{lem}
    \begin{pf}
        Let $\{\varepsilon_n\}_{n = 1}^{\infty}$ be an arbitrary sequence such that $\lim_{n \to \infty} \varepsilon_n = 0$. It is sufficient to prove that $\lim_{n \to \infty} \boldsymbol{\xi}(\varepsilon_n) \to \mathbf{x}^*$. Suppose a contrary, then there exists a neighbourhood $\mathcal{O}$ of $\mathbf{x}^*$ such that the sequence $\{\boldsymbol{\xi}(\varepsilon_n)\}_{n = 1}^{\infty}$ has the infinite number of elements in the compact set $\mathcal{X}\setminus\mathcal{O}$. Hence, there exists a sub-sequence $\{\boldsymbol{\xi}(\varepsilon_{n_k})\}_{k = 1}^{\infty}$ such that
        \begin{equation*}
            \lim_{k \to \infty}\boldsymbol{\xi}(\varepsilon_{n_k}) = \boldsymbol{\xi}^* \in \mathcal{X}\setminus\mathcal{O}.
        \end{equation*}
        On the other hand, we have
        \begin{multline*}
            0 = \lim_{\varepsilon \to 0} f(\boldsymbol{\xi}(\varepsilon)) = \lim_{k \to \infty} f(\boldsymbol{\xi}(\varepsilon_{n_k})) \\
            = f\left(\lim_{k \to \infty}\boldsymbol{\xi}(\varepsilon_{n_k})\right) = f(\boldsymbol{\xi}^*).
        \end{multline*}
        It contradicts to the uniqueness of $\mathbf{x}^*$ as a root of $f(\mathbf{x}) = 0$.
        
    \end{pf}

    \bibliographystyle{plainnat}
    \bibliography{main}
    
\end{document}